%% file: dbNeumann2.tex
\documentclass{article}
\usepackage{latexsym,theorem,graphics,xspace,ifthen}
\usepackage{amssymb,amsmath,amscd}

\include{NMmacros}

\include{MyEnvironments}

\begin{document}
\begin{center}
\large Boundary Regularity for the $\db$-Neumann problem, Part 2 \normalsize\\

\hfill

Robert K. Hladky

\end{center}

\begin{abstract}
We adapt the results of Part 1 to include the unit ball in the Heisenberg group, the model domain with characteristic boundary points. In particular, we construct function spaces on which the Kohn Laplacian with the $\db$-Neumann boundary conditions is an isomorphism. As an application, we establish sharp regularity for a canonical solution to the inhomogenous $\db$ equation on the unit ball.
\end{abstract}

\begin{center}
\tableofcontents
\end{center}

\input{IntroDuction2}

\input{HeiseNberg}
\input{SphericalHarmonics}
\input{GloBal}

\input{BounDary}

\input{INterior}

\input{MaiN}

\appendix
\input{ComParison}

\bibliographystyle{plain}
\bibliography{References}

\end{document}

%% file: NMmacros.tex
\def\today{\ifcase\month\or
 January\or February\or March\or April\or May\or June\or
 July\or August\or September\or October\or November\or December\fi
 \space\number\day, \number\year}

\numberwithin{equation}{section}
\DeclareFontFamily{OT1}{rsfs}{}
\DeclareFontShape{OT1}{rsfs}{m}{n}{ <-7> rsfs5 <7-10> rsfs7 <10-> rsfs10}{} 
\DeclareMathAlphabet{\mathscr}{OT1}{rsfs}{m}{n}

\newcommand{\brckt}[1]{\ifthenelse{\equal{#1}{}}{}{(#1)}}

\newcommand{\del}[2][1]{\ensuremath{\triangle_{[#2]}^{#1}}\xspace}

\newcommand{\tL}{{(\ut,\nu)}}

\newcommand{\db}[1][{}]{\ensuremath{\bar{\partial}^{#1}_{b}}\xspace}
\newcommand{\dbs}[1][{}]{%
\ifthenelse{\equal{#1}{}}{\ensuremath{\bar{\partial}^*_{b}}\xspace}
{\ensuremath{(\bar{\partial}^{#1}_b)^*}\xspace}}%

\newcommand{\dist}[3][U]{\ensuremath{\text{dist}_{#1}\left(#2,#3\right)}\xspace}

\newcommand{\dD}[1][{}]{{\partial_{#1} D}}

\newcommand{\ddDn}[1][\cdot]{\dist[]{\cdot}{\partial D}}

\newcommand{\doD}[1]{\dist{#1}{\dO}}

\newcommand{\dE}[1][{}]{\ensuremath{{d^{#1}_{\mathcal{E}}}}\xspace}

\newcommand{\ndot}{\ensuremath{\!\cdot\!}\xspace}
\newcommand{\Tn}[1]{\ensuremath{T_{\scriptscriptstyle \!#1}}\xspace}
\newcommand{\Tt}{\Tn{\theta}}

\newcommand{\TT}{\Tn{\Theta}}

\newcommand{\Nt}{\ensuremath{N_{\ut}}\xspace}
\newcommand{\NT}{\ensuremath{N_{ \uT}}\xspace}
\newcommand{\Nm}{\ensuremath{N_{\mu}}\xspace}

\newcommand{\cplx}[2]{\ensuremath{\mathcal{B}^{{#1},{#2}}}\xspace}

\newcommand{\norm}[3][{}]{\ensuremath{\big\|#2\big\|_{#3}^{#1}}\xspace}
\newcommand{\snorm}[3][{}]{\ensuremath{\big|#2\big|_{#3}^{#1}}\xspace}

\newcommand{\crT}[1][{}]{\ensuremath{{^o T_{#1}^{\prime}}}\xspace}

\newcommand{\dom}[1]{\ensuremath{\text{Dom}(#1)}\xspace}

\newcommand{\Ker}[1]{\text{Ker}(#1)}

\newcommand{\imag}[1]{\text{Im}\left(#1 \right)}

\newcommand{\be}{{\bar{\eta}}}

\newcommand{\nW}{\ensuremath{\overline{W}}\xspace}

\newcommand{\bt}[1]{{\bar{#1}}}
\newcommand{\Bt}[1]{{\overline{#1}}}

\newcommand{\nY}{\bt{Y}}

\newcommand{\ua}{{\alpha}}
\newcommand{\ub}{{\beta}}
\newcommand{\uc}{{\gamma}}

\newcommand{\ud}{{\delta}}
\newcommand{\ul}{{\lambda}}

\newcommand{\ut}{{\theta}}
\newcommand{\uT}{{\Theta}}
\newcommand{\sT}{{\Theta}}

\newcommand{\vs}{\ensuremath{{\varsigma}}\xspace}

\newcommand{\e}{{\epsilon}}

\newcommand{\dO}{\ensuremath{{\partial \Omega}}\xspace}

\newcommand{\Pt}[1][{\s}]{\ensuremath{{P}^\top_{#1}}\xspace}
\newcommand{\Pb}[1][{\s}]{\ensuremath{{P}^\bot_{#1}}\xspace}

\newcommand{\mcB}[1][{}]{\ensuremath{\mathcal{B}_{#1}}\xspace}

\newcommand{\chs}[1][{w}]{\ensuremath{\mathcal{E}_{#1}}\xspace}
\newcommand{\cho}[1][\Omega]{\ensuremath{\mathcal{E}_{#1}}\xspace}

\newcommand{\bQ}[2][{}]{%
\ifthenelse{\equal{#1}{}}{\ensuremath{\mathcal{Q}_{#2}}\xspace}
{\ensuremath{\mathcal{Q}_{#2}\big( #1 \big)}\xspace}}%

\newcommand{\bI}[1][{}]{%
\ifthenelse{\equal{#1}{}}{\ensuremath{\mathcal{I}}\xspace}{\ensuremath{\mathcal{I}\big(#1\big)}\xspace}}%

\newcommand{\btQ}[2][{}]{%
\ifthenelse{\equal{#1}{}}{\ensuremath{\Bt{\mathcal{Q}}_{#2}}\xspace}
{\ensuremath{\Bt{\mathcal{Q}}_{#2}\big( #1 \big)}\xspace}}%

\newcommand{\CpD}{\ensuremath{C^\infty_+(\overline{D})}\xspace}

\newcommand{\Cic}[1]{\ensuremath{C^\infty_0\brckt{#1}}\xspace}

\newcommand{\CpO}[1][q]{\ensuremath{C^{(0,#1)}_+(\overline{\Omega})}\xspace}

\newcommand{\Ocp}[1][{\gamma}]{\ensuremath{\Omega^+_{#1}}\xspace}
\newcommand{\Ocm}[1][{\gamma}]{\ensuremath{\Omega^-_{#1}}\xspace}
\newcommand{\Ocbp}[1][{\gamma}]{\ensuremath{\Omega^+_{#1}}\xspace}
\newcommand{\Ocbm}[1][{\gamma}]{\ensuremath{\Omega^-_{#1}}\xspace}

\newcommand{\ip}[3]{\ensuremath{\big( \, {#1} \, , \, {#2} \, \big)_{#3}}\xspace}
\newcommand{\ips}[3]{\ensuremath{\big( \, {#2} \, , \, {#1} \, \big)_{#3}}\xspace}
\newcommand{\aip}[3]{\ensuremath{\langle \, {#1} \, , \, {#2} \, \rangle_{#3}}\xspace}

\newcommand{\boxb}[1][{}]{\ensuremath{\square^{#1}_b}\xspace}

\newcommand{\upp}{\ensuremath{^{\prime}}\xspace}
\newcommand{\dupp}{\ensuremath{^{\prime \prime}}\xspace}

\newcommand{\rn}[1]{\ensuremath{\mathbb{R}^{#1}}\xspace}

\newcommand{\sn}[1]{\ensuremath{{\mathbb{S}^{#1}}}\xspace}
\newcommand{\cn}[1]{\ensuremath{\mathbb{C}^{ #1}}\xspace}
\newcommand{\hn}[1]{\ensuremath{\mathbb{H}^{ #1}}\xspace}
\newcommand{\hy}[1]{\ensuremath{\mathbb{U}^{ }}\xspace}

\newcommand{\LMain}[2]{\ensuremath{L^2_{#1}#2}\xspace}

\newcommand{\tlO}{\LMain{\uT}{(\Omega)}}
\newcommand{\lO}[1][\ut]{\LMain{#1}{(\Omega)}}
\newcommand{\lOn}[1][\ut]{\LMain{#1}{}}
\newcommand{\lOa}[1]{\LMain{\tLa}{}}
\newcommand{\tlOn}[1]{\LMain{\uT}{(#1)}}

\newcommand{\Ball}[3]{\ensuremath{B^{#1}_{#2}#3}\xspace}
\newcommand{\ball}[2]{\Ball{#1}{#2}{}}

\newcommand{\bdll}[1]{\ensuremath{D_{\!p}^{#1}}\xspace}

\newcommand{\Lbo}[2][\ut]{\ensuremath{L_{#1}^2\brckt{\tube{#2}}}\xspace}
\newcommand{\Lbp}[3][\ut]{\ensuremath{L_{#1}^2\brckt{\tube[#3]{#2}}}\xspace}

\newcommand{\NewNabMain}[2]{\ensuremath{\nabla^{#1}_{#2}}\xspace}
\newcommand{\NabTotal}[2]{\NewNabMain{#2}{[#1]}}
\newcommand{\Nabx}[2][{}]{\NabTotal{#1}{#2}}
\newcommand{\Nab}[2][{}]{\NewNabMain{#1}{#2}}

\newcommand{\NabT}[1]{\Nab[\sT]{#1}}

\newcommand{\prtl}[2]{\ensuremath{\frac{\partial {#1}}{\partial {#2}}}\xspace}

\newcommand{\hnp}{\ensuremath{\hn{2n+1}\backslash \chs}\xspace}

\renewcommand{\theenumi}{\alph{enumi}}

\newcommand{\WoX}[3]{\ensuremath{\mathring{\mathscr{W}}_{#3}^{#1}#2}\xspace}
\newcommand{\Wo}[1][1]{\ifthenelse{\equal{#1}{1}}%
{\WoX{#1}{(D)}{}}%
{\WoX{#1}{(D)}{\s}}}%

\newcommand{\ED}[1]{\ensuremath{E\brckt{#1}}\xspace}
\newcommand{\V}[1][q]{\ensuremath{\mathscr{V}^{#1}}\xspace}
\newcommand{\Vv}[1][q]{\ensuremath{\mathscr{V}_{\nu}^{#1}}\xspace}
\newcommand{\Vl}[2]{\ensuremath{\mathscr{V}^{#1}_{#2}}\xspace}
\newcommand{\Vvl}[2]{\ensuremath{\mathscr{V}^{#1}_{\nu \Acomma{#2}}}\xspace}

\newcommand{\s}{{\sigma}}

\newcommand{\Gm}[1][\s]{\Gamma\brckt{#1}}

\newcommand{\G}[1][\s]{G(#1)}
\newcommand{\uls}[1][\s]{\ul(#1)}
\newcommand{\Ws}{W_\s}

\newcommand{\nWs}[1][\s]{\Bt{W}_{\!#1}}

\newcommand{\SjMain}[3]{\ensuremath{\mathscr{S}^{#1}_{#2}#3}\xspace}
\newcommand{\SjO}[2][\ut]{\SjMain{#2}{#1}{(\Omega)}}
\newcommand{\Sj}[2][]{\SjMain{#2}{\ut}{\brckt{#1}}}
\newcommand{\Sjn}[2][\ut]{\SjMain{#2}{#1}{}}
\newcommand{\Sjo}[3][\ut]{\SjMain{#2}{#1}{\brckt{#3}}}

\newcommand{\SCo}[3][\ut]{\ensuremath{\mathring{\mathscr{S}}_{#1}^{#2}\brckt{#3}}\xspace}

\newcommand{\Sjb}[4][\theta]{\SjMain{#2}{#1}{\brckt{\tube[#4]{#3}}}}

\newcommand{\HVMain}[3]{\ensuremath{\tilde{\mathscr{S}}_{#1}^{#2}\brckt{#3}}\xspace}
\newcommand{\HV}[2][\ut]{\HVMain{#1}{#2}{\Omega}}
\newcommand{\HVl}[3][\ut]{\HVMain{#1}{#2}{\tube{#3}}}
\newcommand{\HVln}[4][\ut]{\HVMain{#1}{#2}{\tube[#4]{#3}}}
\newcommand{\HVo}[3][\ut]{\HVMain{#1}{#2}{#3}}

\newcommand{\tube}[2][p]{\ensuremath{\Omega_{{#1}}^{#2}}\xspace}

\newcommand{\Acomma}[1]{\ifthenelse{\equal{#1}{}}{}{,#1}}

%% file: MyEnvironments.tex
\newtheorem{thm}{Theorem}[section]

{\theorembodyfont{\rmfamily} }
{\theorembodyfont{\rmfamily} }
{\theorembodyfont{\rmfamily} \newtheorem{rem}[thm]{Remark}}
{\theorembodyfont{\rmfamily} \theoremstyle{plain}  }

\newtheorem{thms}{Theorem}

\newcommand{\pf}{\noindent \textbf{Proof: }}
\newcommand{\epf}{\tiny \ensuremath{\hfill \blacksquare } \normalsize}
\newcommand{\eexm}{\tiny \ensuremath{\hfill \square } \normalsize}

\newcommand{\sName}{none}
\newcommand{\tName}{none}
\newcommand{\eeName}{none}
\newcommand{\eName}{none}

\newcounter{claim}
\newcommand{\setS}[1]{\label{S:#1}\renewcommand{\sName}{#1}}

\newcommand{\bgMain}[5][{}] {%
\renewcommand{\tName}{#4} 
\ifthenelse{\equal{#1}{}}{\begin{#2}\label{#3:#5:\tName}}
{\begin{#2}[#1]\label{#3:#5:\tName}}
}%

\newcommand{\bgT}[2][{}]{\bgMain[#1]{thm}{T}{#2}{\sName}{}}
\newcommand{\enT}{\end{thm}}

\newcommand{\enP}{\end{prop}}

\newcommand{\bgL}[2][{}]{\bgMain[#1]{lemma}{L}{#2}{\sName}}
\newcommand{\enL}{\end{lemma}}

\newcommand{\bgE}[2][{}]{\bgMain[#1]{equation}{E}{#2}{\sName}{}}
\newcommand{\enE}{\end{equation}}

\newcommand{\bgD}[2][{}]{\bgMain[#1]{defn}{D}{#2}{\sName}{}}
\newcommand{\enD}{\end{defn}}

\newcommand{\bgC}[2][{}]{\bgMain[#1]{cor}{C}{#2}{\sName}{}}
\newcommand{\enC}{\end{cor}}

\newcommand{\bgR}[2][{}]{\bgMain[#1]{rem}{R}{#2}{\sName}{}}
\newcommand{\enR}{\end{rem}}

\newcommand{\enA}{\end{appl}}

\newcommand{\enX}{\eexm \end{exm}}

\newcommand{\enQ}{\end{quest}}

\newcommand{\rfT}[2][\sName]{Theorem \ref{T:#1:#2}}
\newcommand{\rfC}[2][\sName]{Corollary \ref{C:#1:#2}}
\newcommand{\rfL}[2][\sName]{Lemma \ref{L:#1:#2}}

\newcommand{\rfR}[2][\sName]{Remark \ref{R:#1:#2}}

\newcommand{\rfS}[1]{Section \ref{S:#1}}

\newcommand{\rfE}[2][\sName]{\eqref{E:#1:#2}}

\newcommand{\bgEn}[2][\tName]{%
\renewcommand{\eName}{#1}
\renewcommand{\theenumi}{#2{enumi}}%
\begin{enumerate}}

\newcommand{\bgEnn}[2][\tName]{%
\renewcommand{\eeName}{#1}
\renewcommand{\theenumii}{$#2{enumii}$}%
\begin{enumerate}}

\newcommand{\enEn}{\end{enumerate}}
\newcommand{\condition}[1]{#1\arabic }

\newcommand{\etem}{\item \label{i:\sName:\eName:\alph{enumi}}}

\newcommand{\rfi}[3][\sName]{(\ref{i:#1:#2:#3})}
\newcommand{\rfI}[1]{\rfi[\sName]{\tName}{#1}}

%% file: IntroDuction2.tex
\section{Introduction}\setS{ID2} 

In Part 1 of this paper, we established existence and sharp regularity for the $\db$-Neumann problem on a specialised class of domains in certain model CR manifolds. The assumption was that the domain could be expressed as the product of a compact normal CR manifold with a precompact open set in the hyperbolic plane. These domains suffered from the severe restriction that they could not possess characteristic boundary points.

In this second part, we shall explore the special case of the unit ball in the Heisenberg group \hn{2n+1}. If the centre of the group is removed, we can exhibit the remaining part of \hn{2n+1} as a product manifold with factors consisting of the unit sphere and the hyperbolic plane. The results of Part 1 can then be employed to establish existence and regularity for the $\db$-Neumann problem on the ball with respect to a class of singularly weighted Folland-Stein spaces. The singularities of the weights occur all along a characteristic axis in \hn{2n+1}. By carefully studying the nature of these spaces and constructing precise interior estimates, we shall be able to establish sharp estimates for solutions to the $\db$-Neumann problem that exhibit singularity only at the characteristic boundary points. 

It is shown in \cite{Folland:S} that the unit sphere has trivial Kohn-Rossi cohomology at the $(0,q)$-level for $1 \leq q \leq n-2$ . Therefore when we apply the results of Part 1, in particular Theorem 9.4,  we see that the Kohn Laplacian is injective and Fredholm.

The main theorem of this paper is as follows: 

\begin{thms}\label{2A} Let $\Omega$ be the unit ball in \hn{2n+1} with $n \geq 4$. Denote by $\dE$ the homogeneous distance of $p$ to the set of characteristic points of the boundary \dO.

\hfill

\noindent Suppose $1 \leq q \leq n-2$. Then on $(0,q)$-forms the operator \[ \boxb\colon  \dE[2] \cdot \HV[\dE]{k,2}\longrightarrow \SjO[\dE]{k}\] is an isomorphism.
\end{thms}
The precise definitions of the spaces involved are very similar in nature to those employed in Part 1 and will be described in detail in \rfS{MN}. This theorem ensures hypoellipticity of \boxb up to all non-characteristic boundary points, but cannot guarantee global hypoellipticity. This is not unexpected as a similar phenomenon occurred in Jerison's study of the Dirichlet problem \cite{Jerison}.

A useful corollary of this theorem is the existence of solutions to the non-homogenous $\db$ equation with sharp estimates. More precisely we shall establish:

\begin{thms}\label{2B} Let $\Omega$ be the unit ball in \hn{2n+1} with $n \geq 4$. Denote by $\dE$ the homogeneous distance of $p$ to the set of characteristic points of the boundary \dO.

\hfill

\noindent Suppose $1 \leq q \leq n-2$. Then for any $(p,q)$-form $\vs$ such that $\db \vs =0$ there exists a $(p,q-1)$-form $\varphi$ such that
\[ \db \varphi = \vs.\]
Furthemore, if $\vs \in \SjO[\sT,\dE]{k}$ then $\varphi$ can be chosen so that $\varphi \in \dE \cdot \HV[\sT,\dE]{k,1}$ and there is a uniform estimate 
\[ \norm{\dE[-1] \varphi}{\HV[\sT,\dE]{k,1}} \leq C \norm{\vs}{\SjO[\sT,\dE]{k}}.\]
\end{thms}

%% file: HeiseNberg.tex
\section{The Heisenberg Group}\setS{HN}

\bgD{Heisenberg}
The Heisenberg group of dimension $2n+1$ is the manifold
\[\hn{2n+1}=\left\{ (t,z) \in \rn{} \times \cn{n} \right\},\] equipped with the CR-structure $\crT$ defined as the complex linear span of the vector fields $L_j =\prtl{}{z^j}+iz^\bt{j} \prtl{}{t}$, $j=1,\dots,n$.
\enD
It is easy to check that the Lie bracket $[L_j,L_k]=0$ for each pair $j,k$. We introduce the function $s=z^k z^\bt{k} = |z|^2$. Then the Heisenberg group can be embedded into \cn{n+1} as a hypersurface by the map $(t,z) \mapsto (w,z)$ where $w=t+is$. The group structure for \hn{2n+1} is induced from the automorphism subgroup of this hypersurface in \cn{n+1}. More concretely
\[ (t_1,z_1) \cdot (t_2,z_2) = (t_1+t_2 +2 \imag{\aip{z_1}{z_2}{}},z_1+z_2).\]

\bgD{Isometries}\hfill
\bgEn{\roman}
\etem The dilation operators of \hn{2n+1} consists of the family $\{ \ud_r \}_{r \in \rn{}_{>0}}$ defined by $\ud_r:(t,z)\mapsto (r^2t,rz)$.
\etem The translation operators $\{\tau_x\}_{x \in \rn{}}$ are defined by $\tau_x:(t,z) \mapsto (t+x,z)$.
\etem $\mathcal{G}$ denotes the group generated by $\{\tau_x\} \cup \{\ud_r\}$.
\etem The group $U(n)$ acts on \hn{2n+1} by $A\cdot(t,z) = (t,Az)$.
\enEn
\enD
It is easy to verify that the elements of $U(n)$ commute with dilations and translations. Furthermore all the elements of $\mathcal{G} \times U(n)$ are CR-diffeomorphisms, i.e. they preserve the bundle \crT.

The function $w=t+is =t+i|z|^2$ introduced earlier is a CR-function. A simple computation shows that the characteristic set of $w$, i.e. the points at which $\db w=0$,  is given by
$ \chs := \{z=0\}.$ Away from \chs the level sets of $w$ are all spheres with a CR structure induced from the inclusion.

The group $\mathcal{G} \times U(n)$ acts transitively on the set $\hnp$. While this is not the full automorphism group, it omits reflections in $t$ for example, it is sufficient for our needs. The orbits of $U(n)$ are the level sets of the CR-function $w$ and $w$ provides a holomorphic coordinate for each orbit of $\mathcal{G}$. 

We define a norm on \hn{2n+1} by 
\bgE{Norm}
\snorm{(t,z)}{H} = \left( |z|^4+t^2 \right)^\frac{1}{4} = \left( s^2 +t^2\right)^\frac{1}{4}.
\enE
This norm is then homogeneous of degree one with respect to the family of dilations.  From this we construct the \textit{homogeneous distance function} on \hn{2n+1} defined by
\[ \dist[H]{(t_1,z_1)}{(t_2,z_2)} = \snorm{(t_1,z_1)\cdot(t_2,z_2)^{-1}}{H} = \left( (t_1-t_2-2\imag{\aip{z_1}{z_2}{}})^2 + |z_1-z_2|^4 \right)^\frac{1}{4}.\]For subsets $K \subset \hn{2n+1}$ we also construct define the distance to $K$ by $\dist[H]{ \cdot }{K} =  \inf\{\dist[H]{\cdot}{p}: p \in K\}$. It is worth noting that if $p \in \chs$ then $[\dist[H]{p}{q}]^2$ is equal to the Euclidean distance from between $w(p)$ and $w(q)$ in $\cn{}$. We shall denote the ball of radius $r$ about a point $p \in \hn{2n+1}$ with respect to the homogeneous distance by $\ball{r}{p}$.

In order to construct an explicit realisation of the $\db$-complex as genuine differential forms and define a Kohn Laplacian, it is necessary to select a pseudohermitian structure for \hn{2n+1}. The choice prevalent throughout the literature is (a constant multiple of) \
\[ \Theta = \frac{1}{2} \left(dt -iz^\bt{k} dz^k + iz^k dz^\bt{k} \right).\] Computation then yields that $d\Theta = idz^k \wedge dz^\bt{k}$ and that the characteristic field is $\TT = 2 \prtl{}{t}$. Since $[L_j,L_\bt{k}]= -i\TT$, this CR-structure is strictly pseudoconvex with the collection $\{L_1,\dots,L_n,L_\bt{1},\dots,L_\bt{n},\TT\}$ forming a global orthonormal frame for $\cn{}TH^{2n+1}$ with respect to the Levi metric \[h_\sT(X,Y) =d\Theta(X,J\bt{Y}) + \Theta(X)\Theta(\bt{Y}).\] The dual coframe is given by $\{dz^1,\dots, dz^n,dz^\bt{1},\dots, dz^\bt{n}, \Theta\}$ and the associated volume form by $dV_\sT=\frac{i^n n!}{2}\, dt \wedge dz^1 \wedge dz^\bt{1} \wedge \dots \wedge dz^n \wedge dz^\bt{n}$. 
 
While the pseudohermitian structure induced by $\Theta$ possesses many exceptional qualities, such as the presence of a global, compatible orthonormal frame, it suffers from the drawback that the orbits of the characteristic field $\TT$  are not closed. Thus it cannot be realised as an example of the structures described in Part 1. 

We can however bring the techniques of Part 1 to bear with the choice of another pseudohermitian form. As has been noted earlier, away from \chs the level sets of the CR function $w$ are all spheres. The induced CR structure on each level set is naturally seen to be CR-diffeomorphic to that on the unit sphere by a simple rescaling. Consider the map $\Xi:\hnp \to \hy{2} \times \sn{2n-1}$, defined by
\[ \Xi: (t,z) \mapsto (t+i|z|^2, \frac{z}{|z|}).\]
Here $\hy{2}$ is the upper half-plane model of hyperbolic 2-space $\{w \in \cn{}: \imag{w} >0\}$ and the manifold on the right is given the smooth product structure. The unit sphere is easily seen to be a normal CR manifold. Therefore we can impose a pseudohermitian structure on $\hy{2} \times \sn{2n-1}$ using the method described in Section 4 of Part 1. This can then be pulled back using the diffeomorphism $\Xi$ to construct a new pseudohermitian structure $\theta$ on $\hnp$. This new pseudohermitian form can be described concretely by
\[\theta=\frac{1}{s}\Theta = \frac{1}{2s} \left(dt -iz^\bt{k} dz^k + iz^k dz^\bt{k} \right).\]
By construction the new characteristic field is tangent to the foliation induced by $w$. It is given explicitly by \[\Tt= i\left(z^k \prtl{}{z^k} - z^\bt{k} \prtl{}{z^\bt{k}}\right).\] In addition $\ut$ is easily seen to be invariant under the action of $\mathcal{G} \oplus U(n)$. Unfortunately while this pseudohermitian structure is now ideal for applying the results of Part 1, it is singular along the characteristic set \chs. For studying domains with purely non-characteristic boundaries this does not matter. Domains with this restriction have been studied in greater generality by R. Diaz \cite{Diaz}. The purpose of this paper is to extend our results to cover the unit ball in \hn{2n+1}. This will only be of real interest if we can obtain information about a Kohn Laplacian that is defined across the characteristic set \chs. We shall therefore examine and compare the Kohn Laplacians associated to both pseudohermitian structures. 

One slight problem with the comparison of the pseudohermitian structures is that they yield different realisations of the abstract tangential Cauchy-Riemann complex $\cplx{p}{q}$ on genuine forms. For a detailed description of the abstract complex , the reader is referred to \cite{Tanaka}. Precisely speaking, a choice of pseudohermitian form induces a global right inverse to the natural projection from $\cn{}\Lambda^{q}$ to $\cplx{0}{q}$. We denote the inverses for our pseudohermitian forms by $\pi^q_\theta$ and $\pi^q_\Theta$ respectively. The composite map $\mu_q = \pi^q_\Theta \circ (\pi^q_\theta)^{-1}$ is then a bijection (away from $\chs$) from realised $(0,q)$-forms for the pseudohermitian form $\theta$ to those determined by $\Theta$. Since the Cauchy-Riemann complex is defined on the abstract quotient bundles,
 the family of operators $\{\mu_q\}$ intertwines $\db[\uT]$ and \db[\ut], i.e. $\db[\uT] \circ \mu_q = \mu_{q+1} \circ \db[\ut]$. In order to establish later regularity results, it will be important to understand the action of $\mu_q$ on $L^2$ spaces and more generally on the Folland-Stein spaces. If we suppose $\Omega$ is a bounded open domain in \hn{2n+1} then for $(0,q)$-forms $\ua$ and $\ub$
\bgE{l}\begin{split}
 \ip{\mu_q \ua}{\mu_q \ub}{\tlO} &= \int\limits_\Omega h_\sT (\mu_q \ua,\mu_q \ub) dV_\sT = \int\limits_\Omega s^{-q} h_\ut(\ua,\ub) s^{n+1} dV_\ut\\
& = \ip{s^{n+1-q} \ua}{\ub}{\lO[\ut]}.
\end{split}\enE

Since the Kohn Laplacian  maps $(0,q)$-forms to $(0,q)$-forms, we can study each value of $q$ independently. We introduce a parameter $\nu \in \rn{}$ and scale the volume form $dV_\ut$ to obtain  $dV_\tL = s^\nu dV_\ut$. When $\nu=n+1-q$ the weighted inner product on $\lO$ matches that of $\lO[\sT]$ on $(0,q)$-forms. Alternatively phrased, when $\nu=n+1-q$ the operator $\mu_q$ is isometric isomorphism from $\lO$ to $\tlO$ on $(0,q)$-forms. Here the space is defined as the $L^2$ space on $\Omega$ corresponding to the volume form $dV_\tL$. When studying non-characteristic domains we shall just consider $\nu$ as a parameter in the range $\nu \geq 0$. For the study of characteristic domains we shall freeze the value of $\nu$ at $n+1-q$ when we study a particular value of $q$. This means that the $L^2$ inner products of the pseudohermitian structures will not agree for forms of degrees away from $q$. In addition since the definition of the Kohn Laplacian involves adjoint operators, it will depend upon the choice of weighting. This has the affect of shifting frequencies in the transverse directions only. This shift will prove beneficial.

The CR equivalence between  \hn{2n+1} and $\hy{} \times \sn{2n-1}$ induces a partial distance function defined by setting $\dist{p}{q}$ to be the hyperbolic distance between $w(p)$ and $w(q)$ in \hy{}.  This concept of distance was used extensively in Part 1 to obtain uniformity of local estimates. Of particular importance were the restricted hyperbolic tubes $\tube{\e}$ defined by
\[ \tube{\e} = \{q \in \overline{\Omega}: \dist{p}{q}<\e\}.\]

For simplicity we shall restrict our attention to the unit ball in \hn{2n+1}, i.e.  we set $\Omega = \{\snorm{(t,z)}{H}<1\}$. Following Part 1, we  identify $\Omega -\chs$ with the product manifold $D \times \sn{2n-1}$ where $D = \{|w|<1\} \subset \hy{2}$. 

We now construct a special function $\varrho$ on $\hn{2n+1}-\chs$. Note that as $D=\{|w|<1\}$, the boundary $\dD$ is a hyperbolic geodesic with the further property that $\dist[]{\cdot}{\dD}$ is preserved by any M\"obius transformation fixing $\dD$. Along the line $t=0$ we can explicitly compute that $\dist[]{is}{\dD} = | \ln s |$. Therefore we can construct a smooth function $\phi(w)$ on \hy{} by setting $\phi(is) = -\ln s$ and declaring that $\phi$ is preserved by all M\"obius transformations that fix $\dD$. 

Now let $\xi$ be a smooth real-valued function on \rn{} such that $\xi(x)=x$ for $|x|<3/4$ and $\xi(x)=1$ for $|x| >1$. On $\hn{2n+1} -\chs$ we define the function $\varrho$ by
\[ \varrho = \xi \circ \phi \circ w.\]

\bgL{Defining}
The function $\varrho$ depends solely upon the real and imaginary parts of $w$ and has the following properties
\bgEn[R]{\condition{R}}
\etem $\varrho$ extends smoothly to $\hn{2n+1}-\cho$ taking the value $1$ on $\Omega \cap \chs$. 
\etem $\varrho=0$ on $\dO$, $\varrho>0$ on $\Omega$ and $\varrho<0$ on $\hn{2n+1}-\overline{\Omega}$.
\etem $\varrho(x) = \doD{x}$ on $\{\doD{x}<3/4\}$ .
\etem $\varrho=1$ on $\{\doD{\cdot} >1\}$.
\etem $b^1 := \inf\limits_{0<\varrho < 3/4}  \min\{|W\varrho|,|\Bt{W}\varrho|\}  >0$.
\etem $B^m := \sup\limits_\Omega  \max\limits_{j+k \leq m} |W^j \bt{W}^k \varrho| < \infty$ for all $m \geq 0$.
\enEn
\enL

\pf
Most of the lemma follows immediately from the construction of $\varrho$. For the bounds on the derivatives we note that the magnitude of $W$ is preserved by M\"obius transformations and so the magnitude of the various derivatives of the function $\dist[]{\cdot}{\dD}$ at any point are determined purely by their values on the line $t=0$. But along $t=0$ we can explicitly compute $\ddDn(t+is)= - \ln s$. By symmetry we this line critical in the $t$-directions and so   $W\ddDn = s \partial_s (- \ln s) =-1$ along $t=0$. An easy induction argument and simple exercise in hyperbolic geometry completes the proof.

\epf

This Lemma shows that the conditions required in Section 6 of Part 1 are satisfied with $\ud(D)=\frac{1}{4}$ and the regularity and estimates derived there apply to this example.

\begin{rem}
The results of this paper will extend to any domain $\Omega$ that admits a function $\varrho$ satisfying the all the above conditions. The modifications for the argument are  slight but add to the technical intricacy. In particular, Theorem \ref{2A} can be extended to any domain $\Omega$ which admits a defining function that depends smoothly on the real and imaginary parts of $w$. The construction of the function $\varrho$ requires only elementary ideas, but is lengthy and adds little insight.
\end{rem}

%% file: SphericalHarmonics.tex
\section{Spherical Harmonics}\setS{SH} 

In \cite{Folland:S} Folland conducted a detailed study of the $\db$-complex on the unit sphere. His results allow for a more refined harmonic analysis  than the frequency decomposition obtained by Tanaka (see Part 1, Section 3 or \cite{Tanaka}) for general normal CR manifolds. Here we present a summary of Folland's results translated into our notation. We recall from Part 1 that for a normal CR manifold there is a decomposition of $L^2$ into joint eigenspaces of $\boxb$ and $-i \nabla_T$. The set of these eigenvectors was denoted  by $\V$ and the respective eigenvalues by $\Gm$ and $\uls$.

\bgT[Folland]{Folland}
For each  $0 \leq q \leq n-1$ there is splitting of each $\V$ into the disjoint union of two pieces $\V[q,\phi]$ and $\V[q,\psi]$, each of which can be decomposed as a disjoint union of finite sets indexed by integers $l,m \geq 0$, and a shift operator \[S:\Vl{q,\phi}{l,m} \to \Vl{q+1,\psi}{l,m-1}, \hspace{1in} (\text{for $0 \leq q \leq n-1$})\] with the following properties: 
\bgEn{\roman}
\etem On $(0,q)$-forms \[L^2(\sn{2n-1}) =  \bigoplus\limits_{\s \in \V} \cn{} \langle \s \rangle.\]
\etem The elements of $\V$ are mutually orthogonal in $L^2(\sn{2n-1})$
\etem For $\s \in \Vl{q,\phi}{l,m}$
\begin{align*}
\db \s&= \begin{cases}  \sqrt{(m+q)(l+n-q-1)} \;  S(\s), \quad &\text{if $q<n-1$;}\\ 0, \quad &\text{otherwise.}\end{cases}\\
\dbs \s &= 0.\\
\Gm &=   (m+q)(l+n-q-1).
\end{align*}
\etem For $\s \in \Vl{q,\psi}{l,m}$, 
\begin{align*}
\db \s &= 0 \\
\dbs \s &= \begin{cases} \sqrt{ (m+q)(l+n-q)}\; S^{-1}(\s), \quad &\text{if $q>0$;}\\
0, \quad&\text{otherwise.}\end{cases}\\
\Gm &= (m+q)(l+n-q)  
\end{align*}
\etem For all $\s \in \Vl{q}{l,m}$, \[ \uls  = l-m-q .\]
\etem For all $l,m$, $\Vl{0,\psi}{l,m} = \Vl{n-1,\phi}{l,m}  =\emptyset.$
\etem For all $l \geq 0$, $\Vl{1,\phi}{l,0} =\emptyset.$
\enEn
\enT
Here we have used the convention that when $\phi$ or $\psi$ is omitted from a set that we are refering to the union of the possibilities, i.e. $\Vl{q}{l,m} = \Vl{q,\phi}{l,m} \cup \Vl{q,\psi}{l,m}$. 

All these results transfer across to the weighted $L^2$ spaces without change. In addition we can explicitly write down each $\Ws$ and $\nWs$, the projections of $\nabla_Y$ and $\nabla_\bt{Y}$ onto each summand of the decomposition. Following the techniques of Part 1 Section 4, we pull these forms back to \hn{2n+1} and add weights appropriate to the relevant $L^2$-structures. This yields a division of the sets \Vv.  

\bgL{W}
For $\s \in \Vvl{q}{l,m}$ we have
\[ \nWs = \nW - \frac{l-m-q+ \nu}{2}, \qquad \Ws = W + \frac{l-m-q -\nu}{2}.\]
\enL
 
Recall from Part 1, that we constructed an explicit constant $\ED{\s}$ such that for each $\s \in \Vv$ 
 \[ \norm[2]{\nWs^* u}{L^2(D)} \geq \ED{\s}\norm[2]{u}{L^2(D)}\]
 for all $u \in \Wo$. This can now be explicitly computed for $\s \in \Vvl{q}{l,m}$ as
 \[ \ED{\s}= \begin{cases} \nu+l-m-q-2, \quad &\text{if $\nu+l-m-q \geq 3$}\\ \frac{(l-m-q+\nu-1)^2}{4}, \quad &\text{if $1<\nu +l -m-q<3$;}\\ \;0, \quad &\text{otherwise.} \end{cases}.\]
From these observations, we can be more precise in determining optimal constants for estimates than was possible in the general theory.
\bgL{Details}
Suppose $\nu \geq 0$ is a fixed constant and $q$ is in the range $1 \leq q \leq n-2$. Then
\bgEn{\arabic}
\etem For any $\s \in \Vv$, $\Gm+\frac{\uls}{2} \geq 0$ and $|\uls| \leq \G[\s]^2$.
\etem  If $\s \in \Vvl{q}{l,m}$ with $q>0$ then $\Gm \geq 2$.
\etem If $\s\in \Vvl{q}{l,m}$ with $q>0$ then $0 \leq \G^2 \leq \min\left\{ 2 \Gm, \, 2\Gm -\uls\right\}$.
\etem If $\tau\in \Vv[0]$ and then $\Gm[\tau]=0$ if $m=0$ and $\Gm[\tau] \geq 2$ otherwise.
\etem If $\tau \in \Vv[0]-\left(\Vvl{0}{0,0} \cup \Vvl{0}{1,0} \right)$ then $\Gm[\tau] +\ED{\tau} \geq \frac{1}{4} $. If $\nu \geq 2$ this holds for all $\tau \in \Vv[0]$. 
\etem If $\tau \in \Vv[0] - \left(\Vvl{0}{0,0}  \cup \Vvl{0}{1,0}\right)$ then $\G[\tau]^2 \lesssim \left\{\Gm[\tau] +\ED{\tau}\right\}$ uniformly.
\etem If $\tau \in \Vv[0]$ and $\nu \geq 2$ then  $\G[\tau]^2 \lesssim \left\{\Gm[\tau] +\ED{\tau}\right\}$ uniformly.
\etem If $\tau \in \Vv[0]-\Vvl{0}{0,0}$ and $\nu \geq 3$ then $ \Gm[\tau] +\ED{\tau} \geq 2$.
\etem If $\tau \in \Vv[0]$ and $\nu \geq 4$ then $\Gm[\tau] +\ED{\tau} \geq 2$
\etem If $\s \in \Vv$ with $q>0$ and $\nu \geq 3$ then $\G \lesssim \left(\Gm - \sqrt{\Gm} -\frac{1}{2}\right)$ uniformly.
\etem If $\tau \in \Vv[0]$ and $\nu \geq 4$ then $\G[\tau] \lesssim \left(\Gm[\tau] -\sqrt{\Gm[\tau]} -\frac{1}{2} + \frac{1}{2}\ED{\tau}\right)$ uniformly.
\etem If $\tau\in \Vv[0]-\Vvl{0}{0,0}$ and $\nu \geq 3$ then $\G[\tau] \lesssim \left(\Gm[\tau] -\sqrt{\Gm[\tau]} -\frac{1}{2} + \frac{1}{2}\ED{\tau}\right)$ uniformly.
\enEn
\enL

\pf
For \rfI{a} we simply write out the expressions in full and the result is obvious.

For \rfI{b} and \rfI{c} we note that $\G^2 \leq 1+\Gm$. It is easy to check that on these summands $\Gm\geq 2$ and $\Gm +2\uls \geq 1$.

 For \rfI{d},\rfI{e} and \rfI{h},\rfI{i} we note that $\Gm[\tau]=m(l+n-1)=0$ if and only if $m=0$. Since $n-1 \geq 2$ we have $\Gm[\tau] \geq 2$ if $m>0$. When $m=0$ then $\ED{\tau} \geq \frac{1}{4}$ when $\frac{\nu +l}{2} \geq 1$ and $\ED{\tau} \geq 2$ when $\frac{\nu+l}{2} \geq 2$.

 For \rfI{f} and \rfI{g} we note that if $m>0$ then it is easy to uniformly bound $\bt{\Gamma}(\tau)\lesssim \Gm[\tau]$ hence get the desired bounds. If $m=0$ then everything is linear in $l$ for large $l$ and uniformly bounded above 0 for small values of $l$.

Finally \rfI{j},\rfI{k} and \rfI{l} all follow easily from previous computations estimating $G$ in terms of $\Gm[\tau]$ and $\ED{\tau}$.

\epf

%% file: GloBal.tex
\section{Global Regularity and Estimates }\setS{GB}

Our goal in this section is to establish local regularity and global basic estimates for the operator $\boxb[\sT]$ on $\Omega$. The main tool are the results of Part 1 applied to the pseudohermitian structure induced by $\ut$. We begin by computing where exactly we can compare the operators \boxb[\ut] and \boxb[\sT]. This amounts to unravelling how their domains are related under the intertwining operators $\mu: \Lambda_\ut \to \Lambda_\sT$. 
\bgL{Domain}
Fix $1 \leq q \leq n-2$ and set $\nu=n+1-q$. Suppose $\vs \in L^2_\sT (\Lambda^{(0,q)}_\sT) \cap  \dom{\boxb[\sT]}$ and $\varphi =\mu^{-1} \vs$.  Then
\bgEn{\roman}
\etem $\varphi \in \dom{\db[\ut]}$ and $(\sqrt{s})^{-1} \db[\ut] \varphi \in \lO$ .
\etem $\varphi \in \dom{\dbs[\ut]}$ and $(\sqrt{s})^{-1} \dbs[\ut] \varphi = \mu^{-1} \left( \sqrt{s} \,\dbs[\sT] \vs \right) \in \lO$.
\etem $\db[\ut] \varphi \in \dom{\dbs[\ut]}$ and $(\sqrt{s})^{-1} \dbs[\ut] \db[\ut] \varphi = \mu^{-1} \left( \sqrt{s} \, \dbs[\sT] \db[\sT] \vs\right)- (\sqrt{s})^{-1} \ut^0 \vee \mu^{-1} \db[\sT] \vs \in \lO$.
\etem $\dbs[\ut] \varphi \in \dom{\db[\ut]}$ and $(\sqrt{s})^{-1}\db[\ut] \dbs[\ut] \varphi = \sqrt{s} \, \mu^{-1} \db[\sT] \dbs[\sT] \vs + \theta^\bt{0} \wedge \sqrt{s}\, \mu^{-1} \dbs[\sT] \vs\in \lO$.
\enEn
\enL

\pf Recall that with our choice $\nu=n+1-q$, the family of operators $\{\mu_k\}$ intertwine the operators $\db[\ut]$ and $\db[\sT]$, but for degrees away from the fixed value of $q$ the $L^2$ spaces only match up to a scale factor.  We can express this as
\[ \ip{\ua}{\ub}{\lO} = \ip{s^{r} \mu \ua}{\mu \ub}{\tlO}\]
when $\ua$ and $\ub$ are $(0,q+r)$-forms.

To show \rfI{a} it suffices to check the $\lO$ integrability of $(\sqrt{s})^{-1} \db[\ut] \varphi$. But
\[ \norm[2]{(\sqrt{s})^{-1} \db[\ut] \varphi}{\lO} = \ip{s (\sqrt{s})^{-1}  \db[\sT] \vs}{(\sqrt{s})^{-1} \db[\sT]\vs}{\tlO} = \norm[2]{\db[\sT] \vs}{\tlO}.\]
For \rfI{b} we first define \CpO to be the set of $(0,q)$ forms in $\lO$ whose partial Fourier functions are all in $\CpD$, i.e. all smooth functions on $\overline{D}$ that vanish in some neighbourhood of the line $\{s=0\}$. It is shown in the appendix of Part 1 that \CpO is a  dense subset of $\dom{\db[\ut]}$  in the graph norm. For all $\ua \in \CpO[q-1]$
\begin{align*}
\ip{\db[\ut] \ua}{\varphi}{\lO} &= \ip{\db[\sT] \mu \ua}{\vs}{\tlO} = \ip{\mu \ua}{\dbs[\sT] \vs}{\tlO}\\
&= \ip{\ua}{s \mu^{-1} \dbs[\sT] \vs}{\lO}.
\end{align*} 
Both sides are continuous in $\ua$ for the graph norm of $\db[\ut]$. This implies that $\varphi \in \dom{\dbs[\ut]}$ and $\dbs[\ut] \varphi = s \mu^{-1} \dbs[\sT] \vs$. Finally to establish \rfI{b} we note that \[\norm[2]{\dbs[\sT] \vs}{\tlO} = \ip{s \mu^{-1} \dbs[\sT] \vs}{\mu^{-1} \dbs[\sT] \vs}{\lO}=\norm[2]{(\sqrt{s})^{-1} \dbs[\ut] \varphi}{\lO}.\]
We follow a similar line of reasoning to establish \rfI{c}. Namely for $\ua \in \CpO$,
\begin{align*}
\ip{\db[\ut] \ua}{\db[\ut] \varphi}{\lO} &= \ip{s\db[\sT] \mu \ua}{\db[\sT] \vs}{\tlO} = \ip{\db[\sT] \left[s (\mu \ua) \right] - s(\mu \ut^\bt{0}) \wedge \mu \ua}{\db[\sT] \vs}{\tlO}\\
& = \ip{s (\mu \ua)}{\dbs[\sT] \db[\sT] \vs}{\tlO} - \ip{s(\mu \ut^\bt{0}) \wedge \mu \ua}{\db[\sT] \vs}{\tlO}\\
&= \ip{\ua}{s \mu^{-1} \dbs[\sT] \db[\sT] \vs}{\lO} - \ip{\ua}{\ut^0 \vee \db[\ut] \varphi}{\lO}.
\end{align*}
Statement \rfI{d} follows easily from \rfI{a} and \rfI{b}.

\epf

\bgC{Standard}
If $1 \leq q \leq n-2$ and $\nu=n+1-q$ then on $(0,q)$ forms,
\[ \dom{\boxb[\sT]} \subset \mu \dom{\boxb[\ut]}.\]
Furthermore if $\vs \in \dom{\boxb[\sT]}$ and $\varphi = \mu^{-1} \vs$ then $(\sqrt{s})^{-1} \boxb[\ut] \varphi \in \lO$ and 
\[\boxb[\ut] \varphi =\mu^{-1} s\boxb[\sT] \vs - \ut^0 \vee \mu^{-1} \db[\sT] \vs + \theta^\bt{0} \wedge  \mu^{-1} s\dbs[\sT] \vs.\]
\enC

An important consequence of this corollary is that the difference between $s\boxb[\sT]$ and \boxb[\ut] is not only a first order operator, but also takes a form that we shall be able to absorb into our estimates. For ease of reference we set \[\boxb[\mu] = \mu^{-1} \circ \boxb[\Theta] \circ \mu\] with the domain defined by $\dom{\boxb[\mu]} = \mu^{-1} \dom{\boxb[\sT]}$. Then \rfC{Standard} implies that $\dom{\boxb[\mu]} \subset \dom{\boxb[\ut]}$ and 
\bgE{MuUt}
 s\boxb[\mu] =  \boxb[\ut] - \theta^\bt{0} \wedge \dbs[\ut] + \theta^0 \vee \db[\ut]
\enE
This can now be written out as
\bgE{Transfer}
s\boxb[\mu] \varphi = \left\{\Pt[] \varphi^\top + \nabla_{\nY} \varphi^\top  - \db[\top] \varphi^\bot\right\}+ \theta^\bt{0} \wedge \left\{ \Pb[] \varphi^\bot - (\nabla_{\nY})^* \varphi^\bot - \dbs[\top] \varphi^\top \right\}
\enE
where $\db[\top]:=\db[\ut] - \ut^\bt{0} \wedge \nabla_{\nY}$ is the pure tangential component of $\db[\ut]$ described in Section 4 of Part 1. Here and throughout the remainder of this section we shall always presuppose that we have set $\nu=n+1-q$ where $q$ is the degree of the forms under consideration.

\bgL{Split}
If $\varphi \in \dom{\boxb[\ut] }$ then $\varphi^\top$, $\theta^\bt{0} \wedge \varphi^\bot \in \dom{\boxb[\ut] }$.
\enL

\pf
By \rfT[SH]{Folland} the Kohn Laplacian on the spheres has trivial kernel for $1\leq q \leq n-2$. Thus the regularity theory of Part 1 implies that if $\varphi \in \dom{\boxb[\ut] }$ then $\varphi \in \HV{0,2}$. This is easily seen to be sufficient.

\epf

\bgL[Regularity]{Regularity}
Suppose $p \in \overline{\Omega}\backslash \chs$ and $0<\e<\ud<\frac{1}{4}$. If $\varphi \in \dom{\boxb[\mu]}$ and $s\boxb[\mu] \varphi \in \Sj[\tube{\ud}]{k}$ then \[ \varphi \in \HVl{k,2}{\e}\]
\enL

Here the upper bound of $\frac{1}{4}$ comes from the construction of $\varrho$ which yielded a natural choice of the $\frac{1}{4}$ for the constant $\ud(D)$ from Part 1.

\pf
The proof follows from an inductive argument. The inclusion of $\dom{\boxb[\mu]}$ in $\dom{\boxb[\ut]}$ implies that $\varphi \in \HV{0,2}$. Thus the statement is true for $k=0$.  

Suppose $k \geq 1$ and that the lemma holds for $j<k$. Choose $\e\upp$  such that $\e < \e\upp < \ud$. Then $\varphi \in \HVl{k-1,2}{\e\upp}$. Combined with $s\boxb[\mu] \varphi \in\Sj[\tube{\ud}]{k}$ this is sufficient by \rfE{Transfer} to show that $\Pt[] \varphi^\top  \in \Sj[\tube{\ud}]{k}$.  The regularity results of Part 1 then yield that $\varphi^\top \in \HVl{k,2}{\e\upp}$. This implies $\dbs[\top] \varphi^\top \in \Sj[\tube{\e\upp}]{k}$ also. Therefore $\boxb[\ut] \varphi \in \Sj[\tube{\e\upp}]{k}$ and by the main regularity theorem of Part 1, $\varphi \in \HVl{k,2}{\e}$.

\epf

We now begin the intricate task of constructing the existence theory and sharp estimates for $\boxb[\mu]$ and $\boxb[\sT]$. Our first step is to obtain zero and partial first order global estimates. We shall actually get improved estimates with singular weights in these cases. Recall from the discussion in \rfS{HN} that the Heisenberg group possesses a nonisotropic norm $|(t,z)| = (t^2+ |z|^4)^{1/4}$ which induces the homogeneous distance function. 
\bgD{Distance}
For any point $x=(t,z) \in \hn{2n+1}$ define $\dE(x)$ to be the smoothed minimum homogeneous distance from $x$ to a characteristic point of the boundary of $\Omega$, i.e. $\dE(x) = \dist[H]{x}{\cho}$ for $x$ near $\cho$ but the function is smoothed to take the value $\frac{1}{2}$ far from $\cho$.
\enD
The smoothing process is necessary as otherwise derivatives of $\dE$ would be discontinuous along the line $t=0$ equidistant from both characteristic points. We shall use this smoothed distance primarily as a weight on derivatives in our various function spaces. In addition to weighting derivative by derivative to obtain the weighted Folland-Stein spaces $\Sj{k,\dE}$ defined in Section 2 of Part 1, we shall need a second type of weighting.
\bgD{Weight}
For a normed space $X \subset L^2$ and function $\phi$ we define
\[ \phi\ndot X = \{\phi u: u \in X\}\]
with norm $\norm{f}{\phi \cdot X} = \norm{\phi^{-1} f}{X}$. 
\enD
These have already shown up in our work. For example, \rfC{Standard} states that if $\mu \varphi \in \dom{\boxb[\sT]}$ then $ \boxb[\ut] \varphi \in \sqrt{s} \ndot \lO$. We now use these ideas to improve the basic estimates derived in Part 1. Recall from Section 3 of Part 1  that the operator $G$ is defined on a $(0,q)$-form by multiplication by $\G$ on each component of the partial Fourier decomposition. 

\bgL{Improved}
Suppose that $n \geq 4$ and $1 \leq q \leq n-2$. If $\varphi \in \dom{\boxb[\mu]}$ then $G^2\varphi$, $\nabla_{\nY} \varphi^\top$ and $(\nabla_{\nY})^* \varphi^\bot$ are all in $\dE[2]\ndot \lO$ and there is a constant $C>0$ independent of $\varphi$ such that
\bgE{imp}\norm{(\dE[-2]) G^2 \varphi}{\lO} + \norm{(\dE[-2]) \nabla_{\nY} \varphi^\top}{\lO} + \norm{(\dE[-2])(\nabla_{\nY})^* \varphi^\bot}{\lO} \leq C\norm{\boxb[\mu] \varphi}{\lO}.\enE
\enL

\pf 
We define a function $p$ on $D$ by $p=w-c$ for some real constant $c<-1$. Thus $p$ is holomorphic and $r=\sqrt{p\bt{p}}$ is bounded above and below on $D$. Both $p$ and $r$ induce maps on $\Omega$ which we shall also denote with the same letters. Some easy computation then yields that \[W \frac{1}{p} = \frac{-2is}{p^2}\]
with a similar formula holding for $\nWs[] \frac{1}{\bt{p}}$. We shall use the decomposition $r^2 = p \bt{p}$ and the fact that $\dom{\nWs^*}$ is stable under multiplication by suitably bounded anti-holomorphic functions. For $\varphi \in \dom{\boxb[\mu]}$  
\begin{align*}
\ip{\frac{1}{r^2}\varphi}{s\boxb[\mu] \varphi}{\lO} 
&= \ip{\frac{1}{r^2} \varphi^\top}{(\Pt + \nabla_{\nY}) \varphi^\top}{\lO}- \ip{\frac{1}{r^2} \varphi^\top}{\db[\top] \varphi^\bot}{\lO}\\
& \qquad  + \ip{\frac{1}{r^2}\varphi^\bot}{(\Pb- (\nabla_\nY)^*) \varphi^\bot}{\lO} \\
& \qquad - \ip{\frac{1}{r^2} \varphi^\bot}{\dbs[\top] \varphi^\top}{\lO}.
\end{align*}

Now decompose $\varphi = \left(\sum\limits_{\s \in \Vv}  \varphi^\top_\s \sigma \right) + \theta^\bt{0} \wedge \left( \sum\limits_{\tau \in \Vv[q-1]} \varphi^\bot_\tau \tau \right)$. Also recall there is a shift operator $S: \Vv[q,\phi] \to \Vv[q+1,\psi]$ such that $\ip{ \s}{\db[\top] \tau}{} = \ip{ \dbs[\top] \s}{\tau}{}=0$ unless $\tau=S(\s)$. We break the computations down in more manageable chunks.
\begin{align*}
\big(\frac{1}{r^2} \varphi^\top,&(\Pt[] + \nabla_\nY) \varphi^\top\big)_{\lO}\\
&=\sum\limits_{\s \in \Vv} \ips{\frac{1}{\bt{p}} (\Pt \varphi^\top_\s + \nWs \varphi^\top_\s)}{\frac{1}{\bt{p}} \varphi^\top_\s}{L^2(D)} \\
& = \sum\limits_{\s \in \Vv} \ips{\frac{1}{\bt{p}} \Gm\varphi^\top_\s + \nWs^* (\frac{1}{\bt{p}} \nWs \varphi^\top_\s) + \frac{1}{\bt{p}}\nWs \varphi^\top_\s}{\frac{1}{\bt{p}}\varphi^\top_\s}{L^2(D)}\\
& = \sum\limits_{\s \in \Vv} \left\{ \Gm \norm[2]{r^{-1} \varphi^\top_\s}{L^2(D)} + \ips{ \frac{1}{\bt{p}} \nWs \varphi^\top_\s }{\nWs \frac{1}{\bt{p}} \varphi^\top_\s}{L^2(D)} + \ips{\frac{1}{\bt{p}} \nWs \varphi^\top_\s}{\frac{1}{\bt{p}} \varphi^\top_\s}{L^2(D)}\right\}\\
&= \sum\limits_{\s \in \Vv} \left\{\Gm \norm[2]{r^{-1} \varphi^\top_\s}{L^2(D)} + \ips{ \frac{1}{\bt{p}} \nWs \varphi^\top_\s }{\frac{1}{\bt{p}} \nWs \varphi^\top_\s + \frac{2is}{\bt{p}^2} \varphi^\top_\s +\frac{1}{\bt{p}}  \varphi^\top_\s}{L^2(D)}\right\} \\
&= \sum\limits_{\s \in \Vv} \left\{ \Gm \norm[2]{r^{-1} \varphi^\top_\s}{L^2(D)} + \norm[2]{\frac{1}{r} \nWs \varphi^\top_\s}{L^2(D)} + \ip{ \frac{\bt{p}+2is}{\bt{p}^2} \varphi^\top_\s}{ \frac{1}{\bt{p}} \nWs \varphi^\top_\s}{L^2(D)} \right\}\\
& = \sum\limits_{\s \in \Vv} \left\{ \Gm \norm[2]{r^{-1} \varphi^\top_\s}{L^2(D)} + \norm[2]{r^{-1} \nWs \varphi^\top_\s}{L^2(D)} + \ip{ \frac{p}{\bt{p}^2} \varphi^\top_\s}{ \frac{1}{\bt{p}} \nWs \varphi^\top_\s}{L^2(D)} \right\}
\end{align*}
Thus we see that 
\bgE{Pt}
\begin{split}
\text{Re} \ip{\frac{1}{r^2} \varphi^\top}{(P^\top + \nabla_\nY) \varphi^\top}{\lO} \geq  \sum\limits_{\s \in \Vv} \Big\{ &\left(\Gm-\frac{1}{2}\right) \norm[2]{r^{-1} \varphi^\top_\s}{L^2(D)} \\
& \qquad + \frac{1}{2} \norm[2]{r^{-1} \nWs \varphi^\top_\s}{L^2(D)} \Big\}.
\end{split}
\enE
By a virtually identical argument, this time commuting across a $p^{-1}$ term we can establish the inequality
\bgE{Pn}
\begin{split}
\text{Re} \ip{\frac{1}{r^2}\varphi^\bot}{P^\bot \varphi^\bot - (\nabla_\nY)^* \varphi^\bot}{\lO} \geq \sum\limits_{\tau \in \Vv[q-1]} \Big\{& \left(\Gm[\tau] -\frac{1}{2} \right) \norm[2]{r^{-1} \varphi^\bot_\tau}{L^2(D)}\\
& \qquad  + \frac{1}{2} \norm[2]{r^{-1} \nWs[\tau]^* \varphi^\bot_\tau}{L^2(D)}\Big\}.
\end{split}
\enE
The remaining terms can be estimated as follows
\begin{align*}
 \Big| \ip{\frac{1}{r^2} \varphi^\top}{\db[\top] \varphi^\bot}{\lO} &+ \ip{\frac{1}{r^2} \varphi^\bot}{\dbs[\top] \varphi^\top}{\lO} \Big|\\
& \leq  \sum\limits_{\tau \in \Vvl{q-1,\psi}{}} \left| \ip{r^{-1} \varphi^\top_{S^{-1}(\tau)}}{r^{-1} \sqrt{\Gm[\tau]} \varphi^\bot_\tau}{L^2(D)} \right| \\
& \qquad + \sum\limits_{\tau \in \Vvl{q-1,\psi}{}} \left| \ips{r^{-1} \sqrt{\Gm[\tau]} \varphi^\top_{S^{-1}(\tau)}}{r^{-1} \varphi^\bot_\tau}{L^2(D)}  \right| \\
& \leq  \sum\limits_{\tau \in \Vvl{q-1,\psi}{}} 2\sqrt{\Gm[\tau]} \norm{r^{-1} \varphi^\top_{S^{-1}(\tau)}}{L^2(D)} \norm{r^{-1} \varphi^\bot_\tau}{L^2(D)}\\
& \leq  \sum\limits_{\s \in \Vvl{q,\phi}{}} \sqrt{\Gm} \norm[2]{r^{-1} \varphi^\top_\s}{L^2(D)} + \sum\limits_{\tau \in \Vvl{q-1,\psi}{}} \sqrt{\Gm[\tau]} \norm[2]{r^{-1} \varphi^\bot_\tau}{L^2(D)} 
\end{align*}
Combining these three computations we see
\begin{align*}
\text{Re} \ip{r^{-2} \varphi}{s\boxb[\mu] \varphi}{\lO} & \geq \sum\limits_{\s \in \Vv} \left\{ \left(\Gm -\sqrt{\Gm} -\frac{1}{2} \right) \norm[2]{r^{-1} \varphi^\top_\s}{L^2(D)}  + \frac{1}{2}\norm[2]{r^{-1} \nWs \varphi^\top_\s}{L^2(D)}\right\}\\
& \qquad + \sum\limits_{\tau \in \Vv[q-1]} \Big\{ \left(\Gm[\tau] -\sqrt{\Gm[\tau]} -\frac{1}{2} \right)\norm[2]{r^{-1} \varphi^\bot_\tau}{L^2(D)}\\
& \qquad + \frac{1}{2} \norm[2]{r^{-1} \nWs[\tau]^* \varphi^\bot_\tau}{L^2(D)}  \Big\}.
\end{align*}
Now recall from the existence theory of \boxb[\ut] that we can improve the estimates on the transverse summands by using the results of Section 5 in Part 1. We obtain
\begin{align*}
\norm[2]{r^{-1} \nWs[\tau]^* \varphi^\bot_\tau}{L^2(D)} &=\norm[2]{(\bt{p})^{-1} \nWs[\tau]^* \varphi^\bot_\tau}{L^2(D)}\\
&= \norm[2]{\nWs[\tau]^* \big((\bt{p})^{-1} \varphi^\bot_\tau\big)}{L^2(D)}\\
& \geq \ED{\tau} \norm[2]{(\bt{p})^{-1} \varphi^\bot_\tau}{L^2(D)}\\
& = \ED{\tau} \norm[2]{r^{-1} \varphi^\bot_\tau}{L^2(D)}.
\end{align*}
It is shown in \rfL[SH]{Details} that 
\bgE{Est}
\G[\tau]^2 \lesssim \left(\Gm[\tau] -\sqrt{\Gm[\tau]} -\frac{1}{2} + \frac{1}{2}\ED{\tau}\right)
\enE
uniformly over $\Vv[q-1]$ provided either $q>1$ or $\nu \geq 4$. Since $\nu=n+1-q$ this is guaranteed by $n \geq 4$.

Furthermore the same lemma shows that $[\G]^2 \lesssim  \left(\Gm - \sqrt{\Gm} -\frac{1}{2}\right)$ uniformly over $\Vv$. Therefore we have a uniform estimate of the form 
\[ \norm[2]{r^{-1} G \varphi}{\lO} \leq C  \norm{r^{-1} \varphi}{\lO} \norm{\frac{s}{r} \boxb[\mu] \varphi}{\lO}.\]
Furthermore this estimate is also uniform over the choice of $c<-1$. Now it is clear that $\left|\frac{s}{r}\right|$ is bounded above by $1$ on $D$ for all choices of $c$. Thus we have the uniform estimate 
\[ \norm{r^{-1} G^2 \varphi}{\lO} \leq C \norm{\boxb[\mu] \varphi}{\lO}.\]
From the above computations we then also obtain uniform estimates of the form 
\begin{align*}
 \sum \norm[2]{r^{-1} \nWs \varphi^\top_\s}{L^2(D)} + \sum \norm[2]{r^{-1} \nWs[\tau]^* \varphi^\bot_\tau}{L^2(D)}  &\leq 2\text{Re} \ip{r^{-1} \varphi}{\boxb[\mu] \varphi}{\lO}+ \norm[2]{r^{-1} \varphi}{\lO}\\
& \leq C\upp \text{Re} \ip{r^{-1} \varphi}{\boxb[\mu] \varphi}{\lO}\\
& \leq C\dupp \norm[2]{\boxb[\mu] \varphi}{\lO}.
\end{align*}

If we let $c \to -1$, the dominated convergence theorem yields that \rfE{imp} holds with $|w+1|^{-1}$ in place of \dE[-2]. Near the characteristic point $(-1,0)$ however we see that $\dE[2] = |w+1|$. We repeat the procedure for the characteristic point at $(1,0)$ by letting $c$ tend to $1$ from the right. Combining these estimates completes the proof.

\epf

\bgR{NeedGlobal}The global estimates on the transverse first order and tangential second order terms  are necessary to circumvent the difficulties of localising low order estimates for $\Pt$. 
\end{rem}

\bgR{n4}
The requirement that $n \geq 4$ is a little surprising. When $n=3$ the estimates for $\norm{\nWs[\tau]^* \varphi^\bot_\tau}{}$ just fail to provide sufficient positivity when $l=m=0$ and we are unable to establish \rfE{Est}. It is still possible to establish Fredholm theorems and some regularity for $(0,1)$-forms when $n=3$, but as of yet I am unable to to obtain sharp estimates. 
\end{rem}

\bgC[Existence]{Exist}
For $n \geq 4$ and $1 \leq q \leq n-2$, the equation $\boxb[\mu] \varphi = \ua \in \lO$ on $(0,q)$-forms has a unique solution $\varphi \in \dom{\boxb[\mu]}$. 

 Furthermore $\dom{\boxb[\mu]} \subset \dE[2] \ndot \lO$.
\end{cor}

\pf
The operator $\boxb[\mu]$ is closed and self-adjoint on $\lO$ as $\boxb[\sT]$ is closed and self-adjoint  on $\tlO$. The estimate of \rfL{Improved} is more than enough to show that $\boxb[\mu]$ is injective with closed range. Self-adjointness then implies $\boxb[\mu]$ is surjective onto $\lO$.

\epf

This  has the immediate corollary 
\bgC{Exist2}

For $n \geq 4$ and $1 \leq q \leq n-2$, the equation $\boxb[\sT] \varsigma = \ub \in \lO[\sT]$ on $(0,q)$-forms has a unique solution $\varsigma \in \dom{\boxb[\sT]}$. 

Furthermore $\dom{\boxb[\sT]} \subset \dE[2] \ndot \lO[\sT]$.
\end{cor}

Our work so far is now sufficient to establish the existence of the Neumann operators on $\Omega$ for our various Laplacians. 

\bgL{Restate} Suppose $\Omega = \ball{1}{0} \subset \hn{2n+1}$ with $n \geq 4$. Fix $1 \leq q \leq n-2$ and set $\nu=n+1-q$.
 
\hfill

Acting on $(0,q)$-forms each of the three Laplacians \boxb[\ut], \boxb[\mu] and \boxb[\sT] admits a Neumann operator, i.e., a bounded inverse mapping $L^2$ into its domain. These operators, denoted by \Nt, \Nm and \NT respectively, have the following properties:
\bgEn{\alph}
\etem Each operator \Nt, \NT and \Nm is injective and surjective onto the domain of the corresponding Laplacian and continuous in the relevant $L^2$-norms.
\etem \Nt and $\Nm$ both map \lO continuously into $\HV{0,2}$.
\etem \Nm maps \lO continuously into \dE[2] \ndot \lO.
\etem The composition $ \boxb[\ut] \circ \Nm$ maps $\lO$ into $(\sqrt{s}) \ndot \lO$.
\etem If  $p \in \overline{\Omega}\backslash \chs$, $0<\e<\ud<\frac{1}{4}$ and $\varphi \in \lO$ such that the restriction of $\varphi$ to $\tube{\ud}$ lies in $\Sj[\tube{\ud}]{k}$ then $\Nt \varphi$ restricts to an element of $\HVl{k,2}{\e}$. Further there is a constant $C>0$ independent of $\varphi$, such that
\[ \norm{\Nt \varphi}{\HVl{k,2}{\e}} \leq C  \norm{\varphi}{\Sj[\tube{\ud}]{k}} + C \norm{\varphi}{\lO}.\]
\enEn
\enL

%% file: BounDary.tex
\section{Local Estimates at the Boundary}\setS{BD}

Due to the change in dimension of tangential contact vector fields, it is very difficult to analyse the Kohn Laplacian in neighbourhoods of the characteristic points themselves. To avoid this, we adopt the technique of covering the interior of the domain with a countable collection of balls which shrink in size approaching the characteristic points. This unfortunately means we cannot employ a finite cover, so we must work harder to obtain estimates on small balls that vary in some uniform sense as the balls approach the characteristic points.

The cases of balls approaching along the boundary and approaching along the characteristic line will be handled very differently. The former is studied by careful examination of estimates in small hyperbolic tubes, the latter by employing interior estimates and the automorphism group. However both arguments will employ some technical results concerning the geometry of $\Omega$ near its characteristic points.

 We define
\[ \Ocbm = \left\{ x \in \overline{\Omega} - \chs: s(x) <(\sin \uc)  \dE(x) \right\}\]
with $\Ocbp$ defined similarly but with the inequality reversed. Thus $\Ocbm$ is the intersection of two cones, each based at a characteristic boundary point of $\Omega$. We fix once and for all a choice of $\uc$ that ensures $\varrho =1$ on $\Ocbm$. This implies that any point $x$ with $\doD{x} < \frac{3}{4}$ must lie in $\Ocbp$.

Recall that we have defined restricted hyperbolic tubes by $ \tube{r}:= \{x \in \Omega: \text{dist}_{U} (p,x)<r\}$. These restricted tubes are then the preimages  of the restricted hyperbolic balls $\bdll{r}$ by the projection from $\Omega$ onto $D$. 

\bgL{Cover}
For all $0<\e<\ud<\frac{1}{4}$ there exists a countable collection $\{p_m\}$ of points in $\Ocbp$ with the following properties:
\bgEn{\arabic}
\etem The restricted tubes $\tube[p_m]{\e}$ cover $\Ocp$.
\etem The collection $\{\tube[p_m]{\ud}\}$ is uniformly locally finite.
\etem There exists a constant $C$ such that for all $m$ 
\[\sup\limits_{\tube[p_m]{\ud}} \left( \dE[2] \left[\inf\limits_{\tube[p_m]{\e}} \dE \right]^{-2} \right) < C.\]
\enEn
\end{lemma}

\pf
Zorn's Lemma implies the existence of maximal collections of points in $ \Ocbm$ such that restricted hyperbolic tubes of radius $\e/2$ about these points are disjoint. Let $\{p_m\}$ be any such maximal collection. Then every $p \in \Ocp$ is contained in some $\tube[p_m]{\e}$. Now fix one such $p_j$ and suppose $\tube[p_k]{\ud} \cap \tube[p_j]{\ud} \ne \emptyset$ for some $k$. Then $\tube[p_k]{\e/2} \subset \tube[p_j]{2\ud+\e}$. The volume of the projection to \hy{2} of an unrestricted hyperbolic tube depends solely on its radius. Since the tubes of radius $\e/2$ are disjoint, it follows that only finitely many can be contained in the tube $\tube[p_j]{2\ud+\e}$. Thus only finitely many of the tubes of radius $\ud$ can intersect $\tube[p_j]{\ud}$. This establishes the existence of a cover with properties \rfI{a} and \rfI{b}.

Property \rfI{c} is just a characteristic feature of these hyperbolic tubes. Fix a characteristic point $p$ and define a function on \hn{2n+1} by 
\[F(x)=\sup\limits_{y \in \tube[x]{\ud}} \left( [\text{dist}_H(p,y)]^2 \left[\inf\limits_{z\in\tube[x]{\e}} \text{dist}_H(p,z) \right]^{-2} \right).\] 
Since $\Omega$ is compact $F(x)$ is bounded above on the set $\{x \in \Ocbm[\ua]: \text{dist}_H(x,p)=1\}$ for any $0<\ua<\pi/2$. Any dilation centred at $p$ on \hn{2n+1} projects down to $\hy{2}$ as an element of the M\"obius group. It therefore must preserve the hyperbolic distance while scaling the homogeneous distance by a scalar factor. It follows that $F$ is homogeneous of degree $0$ under any such dilation. Therefore $F$ is bounded everywhere on $\Omega-\chs$.
 
 \epf

\bgL{Convert}
Choose $0<\e<\ud<\frac{1}{4}$ . Choose any $p \in \Ocbp$. Then there is a constant $C$ depending solely on $\e$, $k$, $\gamma$ and $\ud$ such that whenever $\varphi \in \dom{\boxb[\mu]}$ with $\boxb[\mu] \varphi \in \Sjb{k}{\ud}{p}$, the following estimate holds
\bgE{Convert}
 \begin{split}
\norm{\boxb[\ut] \varphi}{\Sjb{k}{\e}{p}} &\leq C \Big\{ \norm{s \boxb[\mu] \varphi}{\Sjb{k}{\ud}{p}} + \norm{ [G]^2 \varphi}{\Lbo{\ud}} \\
& \qquad + \norm{\nabla_{\bar{Y}} \varphi^\top}{\Lbo{\ud}} + \norm{ (\nabla_{\bar{Y}})^* \varphi^\bot}{\Lbo{\ud}} \Big\} .
\end{split}
\enE
\enL

\pf 
For convenience of reference we denote the sum of the last three terms of the right hand side of \rfE{Convert} by $\mcB[\ud](p)$. The case $k=0$ follows trivially from the definitions. When $k=1$, we first note that if $\e < \ub <\ud$ then
\begin{align*}
\norm{ \boxb[\ut] \varphi^\top}{\Sjb{1}{\e}{p}} &\leq \norm{s\boxb[\mu] \varphi}{\Sjb{1}{\e}{p}} + \norm{\nabla_{\bar{Y}} \varphi^\top}{\Sjb{1}{\e}{p}} + \norm{\dbs[\top] \varphi^\bot}{\Sjb{1}{\e}{p}}\\
& \leq C \norm{s\boxb[\mu] \varphi}{\Sjb{1}{\e}{p}} + \norm{ \boxb[\ut] \varphi^\top}{\Lbo{\ub}} + \mcB[\ub](p)\\
& \leq C\upp \norm{s\boxb[\mu] \varphi}{\Sjb{1}{\ud}{p}} + \mcB[\ud](p)
\end{align*}
by \rfE[GB]{Transfer} and the regularity results of \rfL[GB]{Regularity} and Part 1. 
\begin{align*}
\norm{ \boxb[\ut] \varphi^\bot}{\Sjb{1}{\e}{p}} &\leq \norm{s\boxb[\mu] \varphi}{\Sjb{1}{\e}{p}} + \norm{(\nabla_{\bar{Y}})^* \varphi^\bot}{\Sjb{1}{\e}{p}} + \norm{\db[\top] \varphi^\top}{\Sjb{1}{\e}{p}}\\
& \leq C \norm{s\boxb[\mu] \varphi}{\Sjb{1}{\e}{p}} +C \norm{\boxb[\ut] \varphi^\bot}{\Lbo{\ub}} + C \norm{\varphi^\top }{\Sjb{1}{\ub}{p}}\\
& \leq C\upp \big\{  \norm{s\boxb[\mu] \varphi}{\Sjb{1}{\ud}{p}} + \mcB[\ud](p) \big\} .
\end{align*} 
This establishes the case $k=1$. The remaining cases with $k>1$ are reduced to this one by a similar induction argument. From the local regularity result of Part 1, Theorem 9.11
\begin{align*}
\norm{ \boxb[\ut] \varphi^\top}{\Sjb{k}{\e}{p}} &\leq \norm{s\boxb[\mu] \varphi}{\Sjb{k}{\e}{p}} + \norm{\nabla_{\bar{Y}} \varphi^\top}{\Sjb{k}{\e}{p}} + \norm{\dbs[\top] \varphi^\bot}{\Sjb{k}{\e}{p}}\\
& \leq C \norm{s\boxb[\mu] \varphi}{\Sjb{k}{\e}{p}} + C \norm{\varphi}{\HVl{k-1,2}{\e}}\\
& \leq C\upp \norm{s\boxb[\mu] \varphi}{\Sjb{k}{\e}{p}} +C\upp \norm{\boxb[\ut] \varphi}{\Sjb{k-1}{\ub}{p}} 
\end{align*}
for any $\ub$ with $\e < \ub < \ud$ with the bounding constants depending solely on $\e$ and $\ub$. So by induction on a carefully nested family of open balls we see \[\norm{\boxb[\ut] \varphi^\top}{\Sjb{k}{\e}{p}} \leq C \norm{\boxb[\mu] \varphi}{\Sjb{k}{\ud}{p}} + C \mcB[\ud](p)\]
with the bounding constant depending solely on $\e$ and $k$. Then
\begin{align*}
\norm{ \boxb[\ut] \varphi^\bot}{\Sjb{k}{\e}{p}} &\leq \norm{s\boxb[\mu] \varphi}{\Sjb{k}{\e}{p}} + \norm{(\nabla_{\bar{Y}})^* \varphi^\bot}{\Sjb{k}{\e}{p}} + \norm{\db[\top] \varphi^\top}{\Sjb{k}{\e}{p}}\\
& \leq C \norm{s\boxb[\mu] \varphi}{\Sjb{k}{\e}{p}} +C \norm{\varphi}{\HVl{k-1,2}{\e}} + C\norm{\varphi^\top }{\HVl{k,2}{\e}}\\
& \leq C\upp  \norm{s\boxb[\mu] \varphi}{\Sjb{k}{\e}{p}} + C\upp\norm{\boxb[\ut] \varphi}{\Sjb{k-1}{\ub}{p}} + C\upp\norm{\boxb[\ut] \varphi^\top}{\Sjb{k}{\ub}{p}}
\end{align*} 
for some $\ub$ with $\e<\ub<\ud$. The estimate then follows from the previous computation and induction.

\epf

\bgC{Local}
Under the conditions of \rfL{Convert} , there is an estimate 
\begin{align*}
\norm{\varphi}{\HVl{k,2}{\e}} &\leq C  \Big\{ \norm{s \boxb[\mu] \varphi}{\Sjb{k}{\ud}{p}} + \norm{ [G]^2 \varphi}{\Lbo{\ud}} \\
& \qquad + \norm{\nabla_{\bar{Y}} \varphi^\top}{\Lbo{\ud}} + \norm{ (\nabla_{\bar{Y}})^* \varphi^\bot}{\Lbo{\ud}} \Big\} 
\end{align*}
with the bounding constant depending solely on $\e, \ud$ and $k$. \enC

\pf
Combine the local estimates for $\boxb[\ut]$ with \rfL{Convert}

\epf

%% file: INterior.tex
\section{Local Estimates along the Characteristic Line.}\setS{IN}

We now attend to the task of studying regularity and estimates along the line \chs. Away from the characteristic boundary points, this line should not be problematic for the $\db$-Neumann problem. Indeed the local results of \cite{Folland:H} imply that we have control over a full complement of Folland-Stein derivatives in the interior of $\Omega$. Therefore instead of using the foliation by spheres which degenerates along \chs, we apply these interior estimates to a sequence of small balls that approach the boundary characteristic points. We construct this sequence careful so we can easily utilise the dilation and translation structure of the Heisenberg group.

\bgL{Cover}
There exists a constant $0<\ua<1$ such that $\Ocm$ can be covered by a countable collection of balls $\ball{\ua r_k}{p_k}$ with the following properties:
\bgEn{\arabic}
\etem Each $p_k \in \Omega \cap \chs$.
\etem Each $\ball{r_k}{p_k} \subset \Ocm[2\uc]$.
\etem There is some $\e>0$ such that $(r_k)^2 = [\dE (p_k)]^2 \cos (2\uc)$ when $\dE(p_k) <\e$ and $r_k$ is uniformly bounded below for $\dE(p_k) \geq \e$.
\etem The collection $\ball{r_k}{p_k}$ is locally uniformly finite.
\enEn
\end{lemma}

\pf
Clearly it is sufficient to establish the existence of such a cover of $\Ocm[\uc] \cap U$ for any neighbourhood $U$ of each characteristic point $p$. As there are only two characteristic points these local covers can easily be combined to yield the desired global cover.

Recall that if $p$ is on the characteristic line then the projection to \cn{} of the homogeneous ball of radius $r$ corresponds to the Euclidean ball of radius $r^2$ about the projection of $p$.
 
Without loss of generality we shall consider characteristic point $p=(-1,0)$ and choose $p_0=(t_0,0)$ to be any point in \chs nearer to $p$ than to $(1,0)$.   Thus $\dE(p_k) = \sqrt{t_0+1}$. Set $(r_0)^2 = (1+t_0) \sin (2\uc)$. 

Next set $x= 1 - \dfrac{r_0^2}{t_0+1}$. Inductively define $t_{k+1} = x(1+t_k)-1$ and $r_{k+1} =  r_k\sqrt{x}$ and set $p_{k+1} =(t_{k+1},0)$. Choose any $\ua > (x+1)^{-1}$. Then the balls $\ball{\ua r_k}{p_k}$ and $\ball{\ua r_{k+1}}{p_{k+1}}$ are not disjoint. Since $t_k \to 0$ the whole collection then covers the part of the characteristic line strictly between $p$ and $p_0$. In fact it is easy to see using elementary properties of similar triangles that they do actually cover the part of $\Ocm$ that intersects some neighbourhood of $p$.

\epf

\bgL[Interior Regularity]{Regularity}

Fix  $r>0$ small enough that $\ball{r}{0}$ is contained strictly in the interior of $\Omega$ and choose $\ua$ with $0<\ua<r$. Suppose that the $(0,q)$-form $\vs \in \dom{\boxb[\sT]}$ satisfies $\boxb[\sT] \vs \in \Sjo[\sT,r]{k}{\ball{r}{0}}$. Then $\vs \in \Sjo[\sT,r]{k+2}{\ball{\ua r}{0}}$ and there is an estimate uniform over $\vs$ and $r$ of the form
\[ \norm{r^{-2} \vs}{\Sjo[\sT,r]{k+2}{\ball{\ua r}{0}}} \lesssim \norm{\boxb[\sT] \vs }{\Sjo[\sT,r]{k}{\ball{r}{0}}} + \norm{r^{-2} \vs}{\tlOn{\ball{r}{0}}}.\]
\enL

\pf
Since these nonisotropic balls are strictly contained in the interior of $\Omega$, this is a statement on the regularity of the formal Kohn Laplacian studied in \cite{Folland:H}. It follows from the work of Folland and Stein in this paper that the regularity statement holds. Furthermore they show the following estimate for the (formal) Kohn Laplacian
\bgE{Formal}
 \norm{\vs}{\Sjo[\sT]{k+2}{\ball{\ua}{0}}} \lesssim \norm{\boxb[\sT] \vs }{\Sjo[\sT]{k}{\ball{1}{0}}} + \norm{\vs}{\tlOn{\ball{1}{0}}}.
\enE

The dilation map $\ud_r$ is an isomorphism of pseudohermitian manifolds from $\left( \hn{2n+1}, r^2 \uT \right)$ to $\left(\hn{2n+1}, \uT \right)$. Multiplication of the pseudohermitian form by a scalar does not affect the Webster-Tanaka connection or the associated $\db$ operator. However it does re-scale the volume form so $L^2$ adjoints are affected. It is easily seen that formally $\boxb[(r^2 \sT)] = r^2 \boxb[\sT]$. Hence we can directly compute
\begin{align*} 
\norm{ \ud_r^* \vs}{\Sjo[\sT]{k}{\ball{\ua}{0}}} &= \norm{ \left(1+ \Nabx[H]{\sT}\right)^{k} \ud_r^* \vs}{\tlOn{\ball{\ua}{0}}}= \norm{ \ud_r^* r^q \left(1+r\Nabx[H]{\sT}\right)^k \vs}{\lOn[(r^2\sT)](\ball{\ua}{0})}\\
&= r^q \norm{\left(1+r\Nabx[H]{\sT}\right)^k \vs}{\lOn[\sT](\ball{r\ua}{0})}= r^q \norm{ \vs}{\Sjo[\sT,r]{k}{\ball{\ua r}{0}}}.
\end{align*}
Thus \rfE{Formal} applied to $\ud_r^* \vs$ yields
\[ r^q \norm{\vs}{\Sjo[\sT,r]{k+2}{\ball{\ua r}{0}}} \lesssim r^{q+2} \norm{\boxb[\sT] \vs }{\Sjo[\sT,r]{k}{\ball{r}{0}}} + r^q \norm{\vs}{\tlOn{\ball{r}{0}}}.\]
This completes the proof.

\epf

\bgR{Translate}
The translation operators preserve the pseudohermitian structure exactly. The previous result therefore holds for balls centred at any point $p \in \chs \cap \Omega$ provided $r$ is chosen sufficiently small.
\end{rem}

\bgL{Local}
Choose $p \in \chs \cap \Omega$ and $0 < \ua \leq 1$. Suppose $r^2 = \dE[2]{p} \cos \uc$ for some $0<\uc<\pi$. Then for each $k$ there are constants $c_1,C_1 >0$ depending only on $\uc$ and $\ua$ such that
\[ c_1 \norm{\vs}{\Sjo[\sT,r]{k}{\ball{\ua r}{p}}} \leq \norm{\vs}{\Sjo[\sT,\dE]{k}{\ball{\ua r}{p}}} \leq C_1 \norm{\vs}{\Sjo[\sT,r]{k}{\ball{\ua r}{p}}}.\]
\enL

\pf
On $\ball{\ua r}{p}$ we immediately see that \[   r \sqrt{\sec \uc -\ua^2} \leq \dE \leq \ua r \sqrt{\sec \uc +\ua^2}.\] Near the characteristic points the function $\dE$ agrees with some $d_x(\cdot)$ so we can apply the arguments of  \rfR[CP]{Homogeneous}. 
\epf

\bgC{Local}
Under the same conditions as \rfL{Local}, suppose that the $(0,q)$-form $\vs \in \dom{\boxb[\sT]}$ satisfies  $\boxb[\sT] \vs \in \Sjo[\sT,\dE]{k}{\ball{r}{p}}$. Then $\vs \in \Sjo[\sT,\dE]{k+2}{\ball{\ua r}{p}}$ and there is an estimate uniform over $\vs$, $r$ and $p$ of the type
\[ \norm{\dE[-2] \vs}{\Sjo[\sT,\dE]{k+2}{\ball{\ua r}{p}}} \lesssim \norm{\boxb[\sT] \vs }{\Sjo[\sT,\dE]{k}{\ball{r}{p}}} + \norm{\dE[-2] \vs}{\tlOn{\ball{r}{p}}}.\]
\end{cor}

\pf
The corollary follows after we make the observation that commuting the $\dE[-2]$ past all the $\dE (\Nabx[H]{\sT})$ terms in $\Sjn[\sT,\dE]{k}$ produces an equivalent norm. See also the remarks following \rfL[CP]{Commute}

\epf

\bgC{Minus}
Suppose the $(0,q)$-form $\vs \in \dom{\boxb[\sT]}$ satisfies $\boxb[\sT] \in \Sjo[\sT,\dE]{k}{\Ocm[2\uc]}$, then $\vs \in \dE[-2] \cdot \Sjo[\sT,\dE]{k+2}{\Ocm}$ and there is a constant $C>0$ independent of $\vs$ such that
\[ \norm{\dE[-2] \vs }{\Sjo[\sT,\dE]{k+2}{\Ocm}} \leq C \norm{\boxb[\sT] \vs}{\Sjo[\sT,\dE]{k}{\Ocm[2\uc]}} + C \norm{\dE[-2] \vs}{L^2_\sT(\Ocm[2\uc])}.\]
\enC

%% file: MaiN.tex
\section{An Isomorphism Theorem}\setS{MN}

Before we state our main theorem, we need a few additional comments. Firstly when we split $(0,q)$-forms into tangential and transverse components, we need to be careful as $\ut^\bt{0}$ is only a $(0,1)$-form for the $\ut$ pseudohermitian structure. When we are working with $\sT$, we shall need to use
\[ \bt{\eta} \colon= \mu^{-1} \ut^\bt{0} = \frac{z^kd\bt{z}^k}{s}.\]
Using this form $\bt{\eta}$ in place of $\ut^\bt{0}$ we can now define the spaces $\HV[\sT,\phi]{k,j}$ in a fashion directly analogous to the definition of $\HV{k,j}$ in Part 1. However it is important to keep careful track of the weighting function $\phi$.

As before we introduce the spaces $\mathscr{V}^{k}_\phi \SjO[\sT,\phi]{m}$ inductively by $\mathscr{V}^{0}_\phi \SjO[\sT,\phi]{m} = \SjO[\sT,\phi]{m}$ and 
\[ \mathscr{V}^{k+1}_\phi \SjO[\sT,\phi]{m} =\left \{ \varphi \in \mathscr{V}^{k}_\phi \SjO[\sT,\phi]{m}: \rho \phi \NabT{[H]} \varphi, \phi  \NabT{[\top]}\varphi \in \mathscr{V}^{k}_\phi \SjO[\sT,\phi]{m} \right\}\] with corresponding inductively defined norms. From this we obtain
\bgD{FinalSpace}
\[ \HV[\sT,\phi]{k,j} := \left\{ \varphi  \in \mathscr{V}^{j}_\phi \SjO[\sT,\phi]{k}: \varphi^\bot \in \SCo[\sT,\phi]{k+j}{\Omega}, (q+ \nabla^\sT_{\bt{Y}}) \varphi^\top \in \SCo[\sT,\phi]{k+j-1}{\Omega}  \right\}\]
with norm 
\[ \norm[2]{\varphi}{\HV[\sT,\phi]{k,j}} := \norm[2]{\varphi}{\mathscr{V}^{j}_\phi\SjO[\sT,\phi]{k}} + \norm[2]{\varphi^\bot}{\SjO[\sT,\phi]{k+j}} + \norm[2]{(q+\nabla^\sT_{\bt{Y}}) \varphi^\top}{\SjO[\sT,\phi]{k+j-1}}.\]
\enD
Here we use $\SCo[\sT,\phi]{k}{\Omega} = \SjO[\sT,\phi]{k} \cap \SCo[\sT,\phi]{1}{\Omega}$ where $ \SCo[\sT,\phi]{1}{\Omega}$ is the closure of $\Cic{}$ in $\SjO[\sT,\phi]{1}$. 

\begin{rem}
The presence of the $q+\nabla^\sT_\bt{Y}$ term as compared to the $\nabla^\ut_\bt{Y}$ term used in Part 1 is a consequence of the subtle interaction between the the intertwining operator $\mu$ and the Webster-Tanaka connection. For the pseudohermitian form $\sT$ the tangential $(0,1)$-forms are spanned by $\ut^\bt{j}=d\bt{z}^j -\bt{z^j} \be$. Now $d\ut^\bt{j} = d\bt{z}^j \wedge \be - \bt{z}^j \eta \wedge \be$ and so $\db[\sT] \ut^\bt{j} = \be \wedge \ut^{\bt{j}}$ but $\nabla^\sT_\bt{Y} \ut^\bt{j} =  -\bt{z}^j \be$. Thus $ \db[\sT] \ut^\bt{j}  = \be \wedge (1+\nabla^\sT_\bt{Y})\ut^\bt{j}$. Whereas from Part 1 we know that $ (\db[\ut] \mu \ut^\bt{j} )^\bot = \nabla^\ut_\bt{Y} \mu \ut^\bt{j}$.
\end{rem}

We can now state and prove the main result of this paper.
     
\bgT[Main Theorem]{Main}

\hfill

\noindent Let $\Omega$ be the unit ball $B^1_0= \{\snorm{(t,z)}{H}< 1\} \subset  \hn{2n+1}$ with $n \geq 4$. Denote by $\dE$ the (smoothed) homogeneous distance of $p$ to the set of characteristic points of the boundary \dO.

\hfill

\noindent Suppose $1 \leq q \leq n-2$ and $k \geq 0$. Then on $(0,q)$-forms the operator \[\boxb[\sT]\colon  \dE[2] \cdot \HV[\sT,\dE]{k,2}\longrightarrow \SjO[\sT,\dE]{k}\] is an isomorphism.
\end{thm}

\pf  Fix $q$ in the range $1 \leq q \leq n-2$ and set $\nu=n+1-q$.
 
First we note that the continuity of $\boxb[\sT]$ between these spaces is clear from the definitions. Injectivity was shown in \rfC[GB]{Exist2}. This same corollary shows that $\boxb[\sT]$ is surjective from its domain onto $\lO[\sT]$. Therefore it is sufficient to show that the Neumann operator $\NT$ maps $\SjO[\sT,\dE]{k}$ continuously into the space $\dE[2] \cdot \HV[\sT,\dE]{k,2}$. 

Fix $\uc$ as in \rfS{BD} and construct the associated cover of $\Ocm$. On $\Ocm$ the function $\varrho$ is bounded above and below, thus we can immediately apply \rfC[IN]{Minus} to see
\bgE{Minus}
 \norm{\NT \vs}{ \dE[2] \cdot \HVo[\sT,\dE]{k,2}{\Ocm}}\approx \norm{ \dE[-2] \NT \vs}{\Sjo[\sT,\dE]{k+2}{\Ocm}} \leq C \norm{\vs}{\SjO[\sT,\dE]{k}} + \norm{ \dE[-2] \NT \vs}{\lO[\sT]}.
\enE

Now choose $0<\uc_1<\uc$ and double cover $\Ocp[\uc_1]$ as in \rfL[BD]{Cover}. It is easy to see that all the larger tubes of the double cover are contained in some \Ocp[\ua] for $0<\ua<\uc_1$. Choose a $(0,q)$-form $\vs \in \SjO[\sT,\dE]{k}$. Then $\vs \in \Sjo[\sT,\dE]{k}{\Ocp[\ua]}$. Thus $\mu^{-1} \vs \in \Sjo{k}{\Ocp[\ua]}$. Set $\varphi = \Nm (\mu^{-1} \vs)$. 

A simple commutation argument (see \rfL[CP]{Commute}) shows that if $p$ is the centre of one of the tubes from the cover then
\[
\norm{\dE[-2] \varphi}{\HVl{k,2}{\e}}  \leq C \big(\inf\limits_{\tube{\e}} \dE \, \big)^{-2} \norm{\varphi}{\HVl{k,2}{\e}}.
\]
Applying \rfL[BD]{Convert} to each point $p_m$ in the constructed cover yields the estimate
\begin{align*}
\norm{\dE[-2]\varphi}{\HVl{k,2}{\e}} & \leq C  \big(\inf\limits_{\tube{\e}} \dE\, \big)^{-2} \left\{ \norm{s (\mu^{-1} \vs)}{\Sjb{k}{\ud}{p}}  + \mcB[\ud](p)\right\}\\
& \leq C \big( \sup\limits_{\tube{\ud}} s \, \big) \big(\inf\limits_{\tube{\e}} \dE\, \big)^{-2} \norm{\mu^{-1} \vs}{\Sjb{k}{\ud}{p}} + C \big(\sup\limits_{\tube{\ud}} \dE\, \big)^{-2} \mcB[\ud](p)\\
& \leq C \Big\{ \norm{ \mu^{-1} \vs}{\Sjb{k}{\ud}{p}}  +  \norm{\dE[-2] [G]^2 \varphi}{\Lbo{\ud}} \\
& \qquad + \norm{\dE[-2] \nabla^\ut_{\bt{Y}} \varphi^\top}{\Lbo{\ud}} + \norm{\dE[-2] (\nabla^\ut_{\bt{Y}})^* \varphi^\bot}{\Lbo{\ud}} \Big\}.
\end{align*}
For the second line we used the final property of the cover constructed in \rfL[BD]{Cover} to uniformly bound the term $\big(\inf\limits_{\tube{\e}} \dE\, \big)^{-2} $ by $ \big(\sup\limits_{\tube{\ud}} \dE\, \big)^{-2} $. For the last line we noted that $\big( \sup\limits_{\tube{\ud}} s \, \big) \big(\inf\limits_{\tube{\e}} \dE\, \big)^{-2}$  can be uniformly bounded.

This  estimate is uniform over $p \in \{p_m\}$. We then apply it to each tube of the cover. Using the local uniform finiteness of the cover we then see
\begin{align*}
\norm{\dE[-2] \varphi}{\HVo{k,2}{\Ocp}} & \leq  \sum\limits_m \norm{\dE[-2] \varphi}{\HVln{k,2}{\e}{p_m}}\\
& \leq C \sum\limits_m \Big\{ \norm{ \mu^{-1} \vs}{\Sjb{k}{\ud}{p_m}}+  \norm{\dE[-2] [G]^2 \varphi}{\Lbp{\ud}{p_m}} \\
& \qquad + \norm{\dE[-2] \nabla^\ut_{\bt{Y}} \varphi^\top}{\Lbp{\ud}{p_m}} + \norm{\dE[-2] (\nabla^\ut_{\bt{Y}})^* \varphi^\bot}{\Lbp{\ud}{p_m}} \Big\}\\
& \leq C \Big\{ \norm{ \mu^{-1} \vs}{\Sjo{k}{\Ocp[\ua]}}  +  \norm{\dE[-2] [G]^2 \varphi}{\lO} \\
& \qquad + \norm{\dE[-2] \nabla^\ut_{\bt{Y}} \varphi^\top}{\lO} + \norm{\dE[-2] (\nabla^\ut_{\bt{Y}})^* \varphi^\bot}{\lO} \Big\}.  
\end{align*} 
Recall that from \rfL[GB]{Improved} we see that the last three terms on the right can be bounded uniformly by $\norm{\mu^{-1} \vs}{\lO}$. Thus we have established that\[ \norm{\dE[-2]\varphi}{\HVo{k,2}{\Ocp[\uc_1]}} \leq C \norm{ \mu^{-1} \vs}{\Sjo{k}{\Ocp[\ua]}} + \norm{\mu^{-1} \vs}{\lO}.\]
When we translate this result over into the Folland-Stein spaces associated to the $\uT$ pseudohermitian form, we notice that since we are working on a positive cone we can replace the $\sqrt{s}$ weights by $\dE$.  Thus we have shown
\bgE{Plus}
 \norm{\dE[-2] \NT \vs}{\HVo[\sT,\dE]{k}{\Ocp[\uc_1]}} \leq C \norm{ \vs}{\Sjo[\sT,\dE]{k}{\Ocp[\ua]}} + C\norm{\dE[-2] \NT \vs}{\lO[\sT]}.
\end{equation}
Combine \rfE{Plus} and \rfE{Minus}, then recall that by \rfL[GB]{Restate},  $\NT$ is continuous from $\tlO$ to $\dE[2] \cdot \tlO$. This completes the proof.

\epf

\bgC{Nonhomogeneous}
With $\Omega$ as in \rfT{Main}, for all $(p,q)$-forms $\vs \in \lO[\sT] $ with $1 \leq q \leq n-2$ such that $\db[\sT] \vs = 0$ the equation
\[ \db[\sT] \varphi = \vs \] is solvable for $\varphi \in \lO[\sT] \cap \left( \Ker{\db[\sT]}\right)^\bot$. Furthermore, if $\vs \in \SjO[\sT,\dE]{k}$ then $\varphi$ can be chosen so that $\varphi \in \dE \cdot \HV[\sT,\dE]{k,1}$ and there is a uniform estimate 
\[ \norm{\dE[-1] \varphi}{\HV[\sT,\dE]{k,1}} \leq C \norm{\vs}{\SjO[\sT,\dE]{k}}.\]
\enC

\pf  When $p=0$, we insist that the solution be orthogonal to the kernal of $\db[\sT]$ and procede exactly as in Theorem 9.5 in Part 1.   We merely need to add the appropriate weights and compute the adjoint $\dbs[\sT]$ using the results of \rfS{GB}. 

When $p>0$ we note that $dz^1,\dots,dz^n$, $dw$ provide a global holomorphic frame for all $(1,0)$-forms. Thus any $(p,q)$-form can be uniquely expressed as a sum of the wedge product of a $(0,q)$-forms with $p$ elements of this frame. The theorem then easily follows from the case $p=0$.

\epf

This result establishes hypoellipticity of solutions only up to non-characteristic boundary points. The weighting by $\dE$ allows for some singularity in the solution even when $\vs$ is globally smooth. It is interesting to note that the argument works when $\vs$ exhibits this same type of singularity at the characteristic points.

%% file: ComParison.tex
\section{Comparison of Folland-Stein Spaces}\setS{CP}

In order to use the regularity theory for \boxb to study that for \boxb[\sT], it is necessary to understand how the Folland-Stein spaces associated with the different pseudohermitian forms $\ut$ and $\uT$ relate. Since the form $\ut$ blows up along the line $\chs=\{z=0\}$, this is analogous to comparing Sobolev spaces for different metrics on unbounded domains.

To facilitate this study we set $\del[k]{H} =(\Nabx[H]{\sT})^k - \mu(\Nabx[H]{\ut})^k \mu^{-1}$ and $\del[k]{\top} = (\Nabx[\top]{\sT} )^k - \mu (\Nabx[\top]{\ut})^k \mu^{-1}.$

\bgL{Delta}
There exists a constant $C$ such that for all sufficiently smooth $(0,q)$-forms $\varsigma$
\[s^{q+k}\snorm[2]{\del[k]{H} \varsigma}{\sT} \leq C  \sum\limits_{j<k} s^{q+j}\snorm[2]{\left(\Nabx[H]{\sT}\right)^j \varsigma}{\sT}\] 
everywhere on $\{s>0\}$. 
\enL

\pf
We work primarily in the coframe $\{dz^j\}$ as this has the huge advantage of being orthonormal and flat with respect to the $\uT$ pseudohermitian form. Since the form $\ut$ is well-behaved away from the line $\{s=0\}$ there is a constant $C$ such that 
\bgE{OnSphere}
\snorm[2]{\del[k]{H} \varsigma}{\sT} \leq C  \sum\limits_{0 \leq j<k} \snorm[2]{\left(\Nabx[H]{\sT}\right)^j \varsigma}{\sT}
\end{equation}
holds on the compact set $\{ (0,z) \in \hn{2n+1}: |z|=1\}$. Now $\ut$ is invariant under the group $\mathcal{G}$ and the form $\uT$ is invariant under the translations and rescales under any dilation by $\ud_r^* \uT =r^2\uT$. The same constant $C$ then works for all points $(t,z)$ with $|z|=1$. The Webster-Tanaka connection is unchanged by scalar multiplication of the pseudohermitian form. We can then establish the lemma at an arbitrary point by pulling back \rfE{OnSphere} by a dilation that maps the point onto the surface described above. Since
\[ r^{2(k+q)} \snorm[2]{\del[k]{H} \varsigma}{\sT} = \snorm[2]{\ud_r^* \del[k]{H} \varsigma}{\sT} = \snorm[2]{\del[k]{H} \ud_r^* \varsigma}{\sT}\]
the results holds as stated.

\epf

\bgC{One}
 For sufficiently smooth $(0,q)$-forms $\varsigma$,
\[ \snorm[2]{(\Nabx[H]{\ut})^k \mu^{-1} \varsigma}{\ut} \leq C \sum\limits_{0 \leq j \leq k}  s^{q+j} \snorm[2]{(\Nabx[H]{\sT})^j  \varsigma}{\sT}.\]
\end{cor}

\pf We compute
\begin{align*}
\snorm[2]{(\Nabx[H]{\ut})^k \mu^{-1} \varsigma}{\ut} &= s^{k+q}  \snorm[2]{\mu (\Nabx[H]{\ut})^k \mu^{-1} \varsigma}{\sT} \leq s^{k+q} \left\{ \snorm[2]{(\Nabx[H]{\sT})^k \varsigma}{\sT} + \snorm[2]{\del[k]{H}  \varsigma}{\sT} \right\}\\
& \leq C \sum\limits_{0 \leq j \leq k}  s^{q+j} \snorm[2]{(\Nabx[H]{\sT})^j \varsigma}{\sT}
\end{align*}
where the last line follows from \rfL{Delta}

\epf

\bgC{Two}
For sufficiently smooth $(0,q)$-forms $\varsigma$,
there is a constant $C_2$ such that
\[ s^{q+k} \snorm[2]{(\Nabx[H]{\sT} )^k \varsigma}{\sT} \leq   C_2 \snorm[2]{\left(1+\Nabx[H]{\ut}\right)^j \mu^{-1} \varsigma}{\ut}.\]
\end{cor}

\pf The proof is by induction with the case $k=0$ obvious. Then
\begin{align*}
s^{q+k} \snorm[2]{(\Nabx[H]{\sT} )^k \varsigma}{\sT} & \leq s^{q+k} \left\{\snorm[2]{\mu (\Nabx[H]{\ut})^k \mu^{-1} \varsigma}{\sT} + \snorm[2]{\del[k]{H} \varsigma}{\sT}\right\}\\
& \leq  \snorm[2]{ (\Nabx[H]{\ut})^k \mu^{-1} \varsigma}{\ut} + C\upp  \sum\limits_{j<k} s^{j+q}\snorm[2]{\Nabx[H]{\sT} \varsigma}{\sT}\\
& \leq   C_2 \snorm[2]{\left(1+\Nabx[H]{\ut}\right)^j \mu^{-1} \varsigma}{\ut}
\end{align*}
for some large constant $C_2$.

\epf

Before moving from these pointwise estimates to results on the Folland-Stein spaces we first establish a couple of technical lemmas which reveal some flexibility in our definitions that will frequently be useful.

\bgL{Commute}
Suppose $\phi$ is a smooth function on $\hn{2n+1}\backslash \chs$ that satisfies the condition that for all $m$ 
\bgE{Commute}
\snorm{ \phi^{m-1} (\Nabx[H]{\sT})^m \phi}{\sT} \leq K(m)
\end{equation}
for some $K(m)$ everywhere on $\hn{2n+1}\backslash \chs$. Then there exist constants $c$ and $C$ such that
\[ c \sum\limits_{j \leq k} \snorm[2]{ \phi^j (\Nabx[H]{\sT})^j \ua}{\sT} \leq \snorm[2]{ \left(1+ \phi \Nabx[H]{\sT}\right)^k \ua}{\sT} \leq C \sum\limits_{j \leq k} \snorm[2]{ \phi^j (\Nabx[H]{\sT})^j \ua}{\sT}\]
everywhere on $\hn{2n+1} \backslash \chs$ for all $k \geq 0$ and all sufficiently smooth contravariant tensors $\ua$.
\end{lemma}

\pf The proof is by induction with the cases $k=0,1$ trivial. Suppose the lemma holds for all $j <k$. Then
\begin{align*}
\snorm[2]{ \left(1+ \phi \Nabx[H]{\sT}\right)^k \ua}{\sT} & \lesssim \sum\limits_{j<k} \snorm[2]{ \phi^j (\Nabx[H]{\sT})^j \ua}{\sT} + \sum\limits_{j<k} \snorm[2]{ \phi^j (\Nabx[H]{\sT})^j \phi \Nabx[H]{\sT} \ua}{\sT}\\
& \lesssim  \sum\limits_{j<k} \snorm[2]{ \phi^j (\Nabx[H]{\sT})^j \ua}{\sT} + \sum\limits_{j<k} \sum\limits_{m \leq j} \snorm[2]{ K(m) \phi^{j-m+1} (\Nabx[H]{\sT})^{j-m+1} \ua}{\sT}\\
&\lesssim  \sum\limits_{j\leq k} \snorm[2]{ \phi^j (\Nabx[H]{\sT})^j \ua}{\sT}
\end{align*}
For the other direction we again employ induction noting
\begin{align*}
\snorm[2]{ \phi^k (\Nabx[H]{\sT})^k \ua}{\sT} &\lesssim \snorm[2]{ \phi^{k-1} (\Nabx[H]{\sT})^{k-1} \phi \Nabx[H]{\sT} \ua }{\sT} + \snorm[2]{\phi^{k-1} [(\Nabx[H]{\sT})^{k-1},\phi] \Nabx[H]{\sT} \ua}{\sT}\\
& \lesssim \snorm[2]{ (1+ \phi \Nabx[H]{\sT})^{k-1} \phi \Nabx[H]{\sT} \ua}{\sT} + \sum\limits_{m\leq k-1} K(k-1-m)\snorm[2]{(\Nabx[H]{\sT})^{m} \ua}{\sT}\\
& \lesssim \snorm[2]{ (1+ \phi \Nabx[H]{\sT})^{k} \ua}{\sT}.
\end{align*}
\epf

\bgR{Homogeneous}
Suppose $\phi$ is smooth on $\hn{2n+1}\backslash \chs$ and is either constant or homogeneous of degree $1$ with respect to any dilation centred at some $p \in \chs$. It is easy to check that $\phi$ then satisfies the condition of the previous lemma. In particular the lemma works for $\phi = \sqrt{s}$ and $\phi= \dist[H]{\cdot}{p}$ for $p \in \cho$.
\end{rem}

\bgR{Commute}
By very similar arguments, commuting powers of $\phi$ across terms of type $(1+ \phi \Nabx[H]{\sT})^k$ produces an equivalent pointwise norm. 
\end{rem}

\bgR{Tangential}
Similar results hold for $\Nabx[\top]{\sT}$ by employing identical arguments. In fact \rfE{Commute} is stronger than the equivalent statement for $\Nabx[\top]{\sT}$.
\end{rem}

\bgC{Homogeneous}
Suppose $\phi_1$ and $\phi_2$ satisfy the conditions of \rfL{Commute}. If $0<\left|\frac{\phi_1}{\phi_2}\right|<K$ on some set $U$ then for all $k \geq 0$ there are constants $c,C>0$  depending only on $k$ and $K$ such that
\[ c\snorm[2]{ \left(1+ \phi_2 \Nabx[H]{\sT}\right)^k \ua}{\sT} \leq \snorm[2]{ \left(1+ \phi_1 \Nabx[H]{\sT}\right)^k \ua}{\sT} \leq C\snorm[2]{ \left(1+ \phi_2 \Nabx[H]{\sT}\right)^k \ua}{\sT}\]
holds for all sufficiently smooth $\ua$ at all points of $U$. 
\end{cor}

\pf
This is obvious when the norms are written out in the equivalent form of \rfL{Commute}.

\epf

We can now translate our earlier pointwise estimates into Folland-Stein estimates. First we recall the nature of the volume forms we are using. Fix a value of $q$ with $1 \leq q \leq n-2$ and set $\nu=n+1-q$. Thus $dV_\ut = s^{-q} dV_\sT$.

The following estimate is then an immediate consequence of \rfC{One} and \rfC{Two}  together with the technical lemma just proved.
\bgC{Comparison}
Suppose $U$ is any bounded open set in \hn{2n+1}. Then there are constants $c$ and $C$ depending only on $k$, $U$ and $q\upp$ such that
\bgE{Comparison}
\begin{split}
 c \norm[2]{s^\frac{q\upp -q}{2} \left(1+ \left(\sqrt{s}  \right) \, \Nabx[H]{\sT} \right)^k \varsigma}{\tlOn{U}} &\leq  \norm[2]{ \left(1+\Nabx[H]{\ut} \right)^k\mu^{-1} \varsigma}{\lOn(U)}\\
& \leq C  \norm[2]{s^\frac{q\upp-q}{2} \left(1+\left(\sqrt{s}\right) \, \Nabx[H]{\sT}\right)^j\varsigma}{\tlOn{U}}
\end{split}
\end{equation}
for all sufficiently smooth $(0,q\upp)$-forms $\varsigma$.
\end{cor}

\bgC{Isomorphism}
With $q$ fixed and $\nu=n+1-q$ we get that the intertwining operator $\mu$ is an isomorphism between  $\Sjn{k}$ and $\Sjn[\sT,\sqrt{s}]{k}$ on $(0,q)$-forms. 
\end{cor}

These arguments can easily be adapted to establish a similar equivalence for the weighted Folland-Stein spaces $\Sjn[\ut,\varrho]{k}$ and $\Sjn[\sT,\varrho \sqrt{s}]{k}$ or for purely tangential components of the derivative.